\documentclass[12pt]{article}

\def\st#1#2{{\mathrel{\mathop{#2}\limits_{#1}}{}\!}}

\def\calx{{\mathcal{X}}}

\def\calu{{\mathcal{U}}}

\def\({\left(}
\def\){\right)}

\def\vsp{\vspace*{1,5mm}\\ }

\def\bk{\bigskip }

\def\sk{\smallskip }

\def\n{\noindent }
\def\dd{\displaystyle}

\def\D{{\Delta}}

\def\barr{\begin{array}}
\def\earr{\end{array}}

\def\bit{\begin{itemize}}

\def\eit{\end{itemize}}

\def\D{{\Delta}}

\usepackage{amsmath, amsthm, amscd, amsfonts, amssymb}
\usepackage[usenames]{color}

\usepackage{mathrsfs}

\numberwithin{equation}{section} 

\newtheorem{theorem}{Theorem}[section]
\newtheorem{proposition}[theorem]{Proposition}

\newtheorem{lemma}[theorem]{Lemma}
\theoremstyle{definition}
\newtheorem{definition}[theorem]{Definition}

\newtheorem{remark}[theorem]{Remark}

\def\S{Schr\"o\-din\-ger}

\def\calf{{\mathcal{F}}}

\def\calx{\mathcal{X}}
\def\calz{\mathcal{Z}}

\def\caly{{\mathcal{Y}}}

\def\bbe{{\mathbb{E}}}
\def\bbr{{\mathbb{R}}}
\def\bbp{{\mathbb{P}}}
\def\bbx{{\mathbb{X}}}

\def\1{^{-1}}
\def\vsp{\vspace*{2mm}\\ }

\def\calf{{\mathcal{F}}}

\def\calx{{\mathcal{X}}}
\def\E{{\mathbb{E}}}
\def\rr{{\mathbb{R}}}
\def\nn{{\mathbb{N}}}
\def\9{{\infty}}

\def\lbb{{\lambda}}
\def\a{{\alpha}}
\def\b{{\beta}}
\def\na{{\nabla}}
\def\g{{\gamma}}
\def\wt{\widetilde}

\def\vf{{\varphi}}
\def\oo{{\omega}}
\def\ooo{{\Omega}}
\def\pp{{\partial}}
\def\D{{\Delta}}
\def\vp{{\varepsilon}}
\def\barr{\begin{array}}
\def\earr{\end{array}}

\def\dd{\displaystyle}
\def\bk{\bigskip }
\def\sk{\smallskip}
\def\n{\noindent }

\def\pas{\mathbb{P}\mbox{-a.s.}}

\def\vsp{\vspace*{2mm}\\ }
\def\ff{\forall }
\def\({\left(}
\def\){\right)}
\def\<{\left<}
\def\>{\right>}

\def\wt{\widetilde}

\def\ol{\overline}
\def\ve{{\varepsilon}}

\begin{document}

\begin{center}
{\Large{\bf Optimal bilinear control of nonlinear stochastic Schr\"odinger equations driven by linear multiplicative noise}}
\bigskip\bk

{\large{\bf Viorel Barbu}}\footnote{Octav Mayer Institute of
Mathematics (Romanian Academy)   and Al.I. Cuza University,
700506, Ia\c si, Romania. },
{\large{\bf Michael R\"ockner}}\footnote{Fakult\"at f\"ur
Mathematik, Universit\"at Bielefeld,  D-33501 Bielefeld, Germany.},
{\large{\bf Deng Zhang}}\footnote{Department of Mathematics,
Shanghai Jiao Tong University, 200240 Shanghai, China. Fakult\"at f\"ur
Mathematik, Universit\"at Bielefeld,  D-33501 Bielefeld, Germany. }
\footnote{E-mail address:
vb41@uaic.ro (V. Barbu), roeckner@math.uni-bielefeld.de (M.
R\"ockner), zhangdeng@amss.ac.cn (D. Zhang)}
\end{center}

\bk\bk\bk

\begin{quote}
\n{\small{\bf Abstract.} Here is investigated the bilinear optimal control problem of quantum mechanical systems with final observation
governed by a stochastic nonlinear \S\ equation perturbed by a linear multiplicative Wiener process.
The existence of an open loop optimal control and first order Lagrange optimality conditions are derived,
via Skorohod's representation theorem, Ekeland's variational principle and the existence for
the linearized dual backward stochastic equation.
Moreover, our approach in particular applies to the deterministic case.}


\end{quote}

\vfill

\section{Introduction}

We consider the controlled stochastic system governed by the nonlinear \S\ equation
\begin{equation}\label{equa-x}
\barr{rcll}
idX(t,\xi)&\!\!=\!\!&\D X(t,\xi)dt+\lbb|X(t,\xi)|^{\a-1}X(t,\xi)dt -i\mu(\xi)X(t,\xi)dt \\
&&    +V_0(\xi)X(t,\xi)dt+\dd\sum\limits^m_{j=1} u_j(t)V_j(\xi)X(t,\xi)dt  \\
&&\dd+iX(t,\xi)dW(t,\xi),\ \ \ t\in(0,T), \ \xi\in\rr^d,\vsp
X(0)&\!\!=\!\!&x\ \mbox{ in }\rr^d. \earr
\end{equation}
Here $\lbb=\pm1,\ \a>1,\ V_j\in W^{1,\9}(\rr^d),\ 0\leq j\leq m,$ are real valued functions, $W$~is the Wiener process,
\begin{equation}\label{e2}
W(t,\xi)=\sum^N_{j=1}\mu_j e_j(\xi)\b_j(t),\ t\ge0,\ \xi\in\rr^d,\end{equation}
and
$$\barr{rcl}
u(t)&=&(u_1(t),...,u_m(t)) \in \bbr^m,\ t\in(0,T),\vsp
\mu(\xi)&=&\dd\frac12\sum^N_{j=1}|\mu_j|^2e^2_j(\xi), \xi\in\rr^d,\ d\ge1,\earr$$
with $\mu_j$ purely imaginary numbers (i.e. $Re \mu_j =0$), $e_j(\xi)$ real-valued functions and $\b_j$ independent real Brow\-nian motions on a probability space $(\Omega,\calf,\mathbb{P})$ with natural filtration $(\calf_t)_{t\ge0}$, $1\le j\le N.$ For simplicity, we assume $N<\9$,
but the arguments in this paper easily extend to the case where $N= \9$.

The physical significance of \eqref{equa-x} is well known. $X=X(t,\xi,\oo),$ $\xi\in\rr^d$, $t\ge0,$ $\oo\in\ooo$, represents the quantum state at time $t$, while the stochastic perturbation $iX\,dW$ represents a stochastic continuous measurement via the pointwise quantum observables $R_j(X)=\mu_je_jX$. The functin $V_0$ describes an external potential.

In the conservative case considered in this paper (i.e. $Re \mu_j=0$, $1\leq j\leq N$), $-i\mu X dt + i X dW$ is indeed the Stratonovitch differential.
It follows by It\^o's formula that $|X(t)|^2_{L^2}=|x|^2_{L^2}$, $\ff t\ge0$. Hence, normalizing the initial state $|x|_{L^2}=1$, we have $|X(t)|_{L^2}=1$, $\ff t\in[0,T],$ and so,
the quantum system evolves on the unit ball of $L^2$ and verifies the conservation of probability. See e.g. \cite{4,5}.

We also mention that, for the general case when $\mu_j$ are complex numbers, one of the main feature is that the mean norm square $|X(t)|_2^2$, $t\in [0,T]$, is a continuous martingale. This fact enables one to define the ``physical'' probability law
and implies the conservation of $\bbe |X(t)|_2^2$, $t\in[0,T]$, which plays an important role in the application to open quantum systems. See, e.g., \cite{bg} for more details.
See also \cite{1,2} for global well-posedness with exponents of the nonlinearity in the optimal subcritical range.

As regards the real valued input control $u$, in most situations it represents an external applied force due to the interaction of the quantum system with an electric field or a laser pulse applied to a quantum system.

Here we shall study an optimal control problem associated with the control system \eqref{equa-x} which, in a few words, can be described as follows (see Problem  (P) below):  find an input control $u$ that steers in time $T$  the state $X$ as close as possible of a target state $\bbx_T$ and a given trajectory $\bbx_1$, and with a reasonable minimum energy.
Roughly speaking, this means to find the quantum mechanical potential $u$ from observation of the quantum state $X(T)$ at the end of time interval $[0,T]$.

It should be mentioned that, there is an extensive literature on the deterministic bilinear control equation \eqref{equa-x} mainly concerned with exact controllability in time $T$ of \S\ equations or with the optimal control problem (see, for instance, \cite{5b}, \cite{5a}, \cite{11a}, \cite{11}, \cite{12},  \cite{13}).
However, there are very few results on optimal control problems governed by nonlinear \S\ equations and,
to the best of our knowledge, none for stochastic control systems  \eqref{equa-x}
with linear multiplicative noise. In the latter case, the existence of an optimal control is largely an open problem,
since the cost functional is not simultaneously lower semicontinuous  and coercive in the basic control space.

The approach we used here is based on Skorohod's representation theorem and Ekeland's variational principle,
and this is one of the main novelties of this work.
The approach is also based on an existence result of the linearized backward dual stochastic equation, which
is also new in the literature and uses sharp stochastic estimates for linear Schr\"odinger equations with time dependent coefficients (see \cite{1,2}).
As a matter of fact, a great effort of this work is dedicated to this issue.

\section{Formulation of problem and the main results} \label{PROBLEM-RESULT}
\setcounter{equation}{0}

To begin with, we recall the definition of a strong solution to equation \eqref{equa-x} (see \cite{1}, \cite{2}).

\begin{definition}\label{d2.1}\rm Let $x\in L^2$ (resp. $H^1$),  $0<T<\9$. Let $\a$ satisfy $1<\a<1+\frac 4 d$ (resp. $1<\a<1+\frac{4}{(d-2)_+}$), $d\geq 1$.
A strong $L^2$-(resp. $H^1$-)solution to \eqref{equa-x} on $[0,T]$ is an $L^2$-(resp. $H^1$-)valued continuous $(\calf_t)_{t\ge0}$-adapted process $X=X(t)$ such that $|X|^{\a-1}X\in L^1(0,T;H^{-1}),$ and $\pas$,
\begin{align}\label{e2.1}
X(t) =&  x-\int^t_0 \bigg( i\Delta X(s)+\mu X(s)+\lbb i|X(s)|^{\a-1} X(s,\xi) +iV_0(\xi)X(s) \nonumber \\
      & \qquad\qquad  +i\dd\sum^m_{j=1}u_j(s)V_j(\xi)X(s) \bigg)ds
          + \int^t_0X(s)dW(s),\  t\in[0,T],
\end{align}
as an It\^o equation in $H^{-2}$ (resp. $H^{-1}$).
\end{definition}
It is easy to check that $\int^t_0X(s)dW(s)$ in Definition \ref{d2.1} is an $L^2$-(resp. $H^1$-)valued continuous stochastic integral.
(We refer, e.g., to \cite{6, LR15} for the general theory of infinite dimensional stochastic integrals.)

Following \cite{1,2}, we introduce the hypotheses below.
\begin{enumerate}
   \item[(H0)] $1<\a<1+ \frac 4 d$. For each $1\leq j\leq N$, $e_j\in C^\9_b(\rr^d)$ satisfies
\begin{align} \label{decay}
   \lim_{|\xi|\to\9}\zeta(\xi)|\pp^\g e_j(\xi)|=0,
\end{align}
where  $\g$ is a multi-index such that $1\leq |\g|\le 2$, and
$$\zeta(\xi)=\left\{\barr{ll}
1+|\xi|^2,&\mbox{if }d\ne2,\vsp
(1+|\xi|^2(\ln(3+|\xi|^2))^2,\ &\mbox{if }d=2.\earr\right.$$

\item[(H1)] In the defocusing case $\lbb=-1$, $1 < \a < 1+ \frac{4}{(d-2)_+}$, and in the focusing case $\lbb=1$,
$1< \a < 1+ \frac{4}{d}$. For each $1\leq j\leq N$, $e_j\in C^\9_b(\rr^d)$
satisfies \eqref{decay} for any multi-index $1\leq |\g| \leq 3$.
\end{enumerate}

The global existence, uniqueness and uniform estimates of the solution to \eqref{equa-x} used in this paper
are summarized in Proposition \ref{Pro-Equ-X} below.

\begin{proposition} \label{Pro-Equ-X}
Assume $(H0)$ (resp. $(H1)$).
For each $x\in L^2$ (resp. $H^1$), $u\in \calu_{ad}$ and $0<T<\9$, there exists a unique strong $L^2$-(resp. $H^1$-)solution $X^u$ to  \eqref{equa-x},
satisfying $|X(t)|_2 = |x|_2$, $t\in [0,T]$ (resp. for any $\rho \geq 1$,
\begin{align} \label{bdd-Xn-H1-Wpq}
     \sup\limits_{u\in \calu_{ad}}\ \bbe \|X^u\|^\rho_{C([0,T]; H^1)}  < \9).
\end{align}

Moreover, assuming  that the exponent $\a$ is in the range specified in $(H1)$ and that $e_k$ are constants, $1\leq k\leq N$, we have for any $\rho \geq 1$,
\begin{align} \label{Cons-bdd-Xn-Lpq}
   \sup\limits_{u\in \calu_{ad}} (  \|X^u\|_{L^\9(\Omega; L^q(0,T; L^{p}))} + \|X^u\|_{L^\rho(\Omega; L^q(0,T; W^{1,p}))}) < \9,
\end{align}
where $(p,q)$ is any Strichartz pair, i.e., $(p,q) \in [2,\9] \times [2,\9], \frac{2}{q} = \frac{d}{2} - \frac{d}{p}$, if $d\not = 2$,
or $(p,q) \in [2,\9) \times (2,\9], \frac{2}{q} = \frac{d}{2} - \frac{d}{p}$, if $d = 2$.
\end{proposition}

The global existence and uniqueness can be proved similarly as in \cite{1,2} by the rescaling approach and the
Strichartz estimates for lower order perturbations of the Laplacian.
We refer to  \cite[Lemma 4.1]{1} and \cite[Lemma 2.7]{2} for explicit formulations of Strichartz estimates in the $L^p$ and Sobolev spaces respectively.
The technical proof
of the estimates \eqref{bdd-Xn-H1-Wpq} and \eqref{Cons-bdd-Xn-Lpq} is postponed to the Appendix for simplicity of the exposition.

In the following,  let
$L^2_{ad}(0,T;\rr^m)$ denote the space of all $(\calf_t)_{t\ge0}$-adapted $\rr^m$-valued processes $u:[0,T]\to\rr^m$ such that $u\in L^2((0,T)\times\ooo;\rr^m)$.
Similarly,  $L^2_{ad}(0,T;L^2(\ooo;L^2))$ denotes the space of $L^2$-valued $(\calf_t)_{t\ge0}$-adapted processes $u$ such that  $\E\int^T_0|u(t)|^2_2dt<\9.$

The optimal control problem we study in the following is
\bit\item[(P)] {\it Minimize $$\E\(\!|X(T)-\mathbb{X}_T|^2_2+\g_1\!\dd\int^T_0\!\!|X(t)-\bbx_1(t)|^2_2dt
+\!\dd\int^T_0\!\! (\g_2|u(t)|^2_m + \g_3|u'(t)|^2_m)  dt\!\)$$\sk\\ on all $(X,u)\in L^2_{ad}(0,T; L^2(\Omega; L^2)) \times \calu_{ad}$\sk\ subject to \eqref{equa-x}.}\eit
Here, $\g_j\ge0$, $1\leq j\leq 3$,  $\bbx_T\in L^2(\ooo,\calf_T, \mathbb{P} ;L^2)$ and $\bbx_1\in L^2_{ad}(0,T;L^2(\ooo;L^2))$
are given. In most situations, $\bbx_1$ is a given trajectory of the uncontrolled system or,  in particular, a steady state solution.
The admissible set $\calu_{ad}$ is defined by
\begin{align}\label{e2.4}
\calu_{ad} =&\bigg\{u\in L^2_{ad}(0,T;\rr^m);\ u\in U,\ \ a.e.\ on\ (0,T)\times\ooo. \bigg\},
\end{align}
where $U$ is a compact  convex subset of $\rr^m$. Let $D_U$ denote the diameter of $U$.
Then, $\sup_{u\in \calu_{ad}}\|u\|_{L^\9(0,T; \bbr^m)} \leq D_U <\9$.

As seen earlier, due to the conservation of $|X(t)|_2^2$, by normalizing the initial state
we have $ |X(t)|_{2}=1$, and so Problem (P) reduces to
\begin{align*}
    &\st{(u,X)}{\rm Min} \bigg\{  - 2 \bbe Re \<X(T),\bbx_T\>_2
              -2 \g_1 \int_0^T Re  \<X(t), \bbx_1(t)\>_2 dt \\
              &\qquad  \qquad \ +  \int_0^T(\g_2 |u(t)|_m^2 + \g_3 |u'(t)|_m^2  ) dt .
    \bigg\}
\end{align*}

It should be said that in the quantum model $V$ is a given potential which describes the spatial profile of an external field,
while the control input $u =\{u_j \}^m_{j=1}$ is its intensity. The objective of the control process is to steer the quantum system from an initial state $x$ to a target state $\bbx_T$ and also in the neighborhood of a given trajectory $\bbx_1$. The last term in the cost functional is the energy cost to obtain the desired objective.

Taking into account that in quantum mechanics the wave function $X$ is not a physical observable,
a more realistic situation is where in the cost functional $|X(T)-\bbx_T|^2_2$ is replaced by $\<Q(X(T))-\bbx_T,X(T)-\bbx_T\>_2$,
where $Q$ is a self-adjoint operator in $L^2$. However, its treatment is essentially the same.

By $\Phi:L^2_{ad}(0,T;\rr^m)\to\rr$ we denote the objective functional
\begin{align} \label{def-Phi}
    \Phi(u)
    =&  \bbe |X^u(T)-\bbx_T|_2^2
      + \g_1 \bbe \int_0^T |X^u(t)-\bbx_1(t)|_2^2 dt   + \g_2 \bbe \int_0^T |u(t)|_m^2dt \nonumber \\
     & +\g_3  \bbe \int_0^T |u'(t)|_m^2 dt,
\end{align}
we may reformulate Problem (P) as
\begin{equation}\label{e2.5}
{\rm(P)}\ \ \ {\rm Min}\{\Phi(u);\ u\in\calu_{ad},\ X^u\ satisfies\  \eqref{equa-x}\}.
\end{equation}

It should be said that, since Problem (P) is a nonconvex minimization problem, in general it is not well posed. However, if $\g_2, \g_3 =0$, we have the following generic existence result.

\begin{proposition}\label{p2.3}
Assume Hypothesis $(H0)$. Then, there is a residual set $$ \mathcal{G}\subset L^2(\ooo,\calf_T,\mathbb{P},L^2)\times L^2_{ad}(0,T;L^2(\ooo;L^2))$$
such that, for every $(\bbx_T,\bbx_1)\in \mathcal{G}$, problem {\rm(P)} has at least one solution $u\in \calu_{ad}$.
\end{proposition}

\n This is an immediate consequence of a well-known result of Edelstein \cite{7} on existence of nearest points of closed sets in uniformly convex Banach spaces. Indeed, if we set $\mathcal{Y}=\{Y=(X^u(T),X^u);\ u\in\calu_{ad}\}$, it follows  that $\mathcal{Y}$ is a closed subset of $L^2(\ooo;\calf_T,\mathbb{P},L^2)\times L^2_{ad}(0,T;L^2(\ooo;L^2))$
(see e.g. the proof of Lemma \ref{Lem-Xn*} and \ref{Lem-conv}) and so,  rewriting Problem {\rm(P)} as
$${\rm Min}\{ \|(\bbx_T,\bbx_1)-Y\|^2_*;\ Y\in \mathcal{Y}\},$$where $\|\cdot\|_*$ is the norm of $L^2(\ooo;\calf_T,\mathbb{P},L^2)\times L^2_{ad}(0,T;L^2(\ooo;L^2)),$ we arrive
at the desired conclusion.
	
However, for the general cases $\g_2, \g_3 \not =0$, the existence of a solution in Problem $(P)$ does not follow by
standard compactness techniques used in deterministic optimization problems. The main reason is that, even if a space $\caly$ is compactly
imbedded into another space $\mathcal{Z}$, one generally does not have the compact imbedding from $L^p(\Omega; \caly)$ to $L^p(\Omega; \mathcal{Z})$, $1\leq p\leq \9$.
Here, we consider the existence for relaxed versions of Problem (P) to be defined below.

\begin{definition}\label{d2.4}
Let $\caly:= L^2(\bbr^d) \times L^2((0,T)\times \bbr^d) \times C([0,T]; \bbr^N) \times L^2(0,T; \bbr^m) \times L^2( (0,T) \times \bbr^d)$
and $(\Omega^*, \mathcal{F}^*, (\mathcal{F}_t^*)_{t\geq 0})$ be a new filtered probability space,
carrying $(\bbx^*_T, \bbx_1^*, \beta^*, u^*, X^*)$ in $\mathcal{Y}$.
Define $L^2_{ad^*}(0,T; L^2(\Omega; L^2))$, $\calu_{ad^*}$ and $\Phi^*(u^*)$  similarly as above on this new filtered probability space.

The system
$(\ooo^*,\calf^*,\mathbb{P}^*,(\calf^*_t)_{t\ge0},\beta^*, u^*,X^*)$ is said to be {\it admissible}, if
$\bbx^*_T \in L^2(\Omega, \mathcal{F}^*_T, \bbp^*; L^2)$, $\bbx^*_1 \in L^2_{ad^*}(0,T; L^2(\Omega; L^2))$,
$\beta^* = (\beta^*_1, \ldots, \beta^*_N)$ is an $(\calf^*_t)_{t\ge0}$-adapted $\bbr^N$-valued Wiener process,
the joint distributions of $(\bbx^*_T, \bbx^*_1, \beta^*)$ and $(\bbx_T, \bbx_1, \beta)$ coincide,
$u^*\in \mathcal{U}_{ad^*}$, and $X^*$ is an $(\calf^*_t)_{t\ge0}$-adapted $L^2$-valued process that satisfies equation \eqref{equa-x}
corresponding to $(\beta^*, u^*)$.

The admissible system $(\ooo^*,\calf^*,\mathbb{P}^*,(\calf^*_t)_{t\ge0},\beta^*, u^*,X^*)$
is said to be a {\it relaxed solution} to the optimal control problem (P),
if
\begin{equation}\label{e2.6}
\Phi^*(u^*)\le\inf\{\Phi(u);\ u\in\calu_{ad},\ X^u\ satisfies\ \eqref{equa-x}\}.
\end{equation}
\end{definition}

We first prove that, under the  regular condition of controls (i.e., $\g_3>0$),
there exists a relaxed solutions for the exponents of the nonlinearity in
exactly the mass-subcritical range.
A similar problem was
studied in \cite{11a} in the deterministic case. We have

\begin{theorem}\label{t2.5}
Consider $\Phi$ with $\g_3>0$.  Assume $(H0)$.
Then, for each $x\in L^2$, $0<T<\9$, there exists at least one
relaxed solution in the sense of Definition \ref{d2.4} to the optimal problem $(P)$.
\end{theorem}

The proof is mainly based on the Skorohod representation theorem and  pathwise analysis
of solutions by the rescaling approach  devoloped in
\cite{1}. We would also like to mention that the rescaling approach allows to obtain pathwise continuous dependence of solutions on controls.\\

In order to construct a relaxed solution with  equality in \eqref{e2.6} in the more difficult irregular
case (i.e., $\g_3=0$), we will employ the Ekeland principle and work with
the dual  backward stochastic equation below
\begin{align}\label{back-equa}
 &d Y= -i\Delta Y\,dt - \lbb i h_1(X^u)Y dt +\lbb i  h_2(X^u)  \ol{Y} dt + \mu Y dt - iV_0Y dt  - i u\cdot V Y dt\nonumber  \\
 &\qquad \quad        + \g_1 (X^u - \bbx_1) dt  - \sum\limits_{k=1}^N \ol{\mu_k} e_k Z_k dt +  \sum\limits_{k=1}^N Z_k d\beta_k(t),    \\
 & Y(T) = -(X^u(T)-\bbx_T),  \nonumber
\end{align}
where Im denotes the imaginary part, and
\begin{equation}\label{e4.4}
h_1(X^u):= \frac{\a+1}{2} |X^u|^{\a-1},\ \ h_2(X^u):=  \frac{\a-1}{2} |X^u|^{\a-3} (X^u)^2 .
\end{equation}
$h_j$, $j=1,2$, are the complex derivatives of the complex function $z\to |z|^{\a-1}z$,
i.e. $h_1(z) = \partial_z (|z|^{\a-1}z)$ and $h_2(z)= \partial_{\ol{z}} (|z|^{\a-1}z)$, $z\in \mathbb{C}$.

However, the singular coefficient $h_{2}(X^u)  $  in
\eqref{back-equa} and the weak regularity effect of the Schr\"odinger group make it quite difficult to obtain the existence
and integrability of the backward
solution. The standard method to derive a global estimate for $\bbe \|Y\|^2_{C([0,T]; L^2)}$
from the It\^o formula applied to $|Y(t)|_2^2$ are not applicable in the nonlinear case.

The idea here is to apply  duality analysis to reduce the
analysis of the backward stochastic equation  to that of the dual equation \eqref{forward-equa-Psi}
below (see also the equation of variation \eqref{forward-equa} below).  By virtue of the forward character of the dual equation,
we will apply the rescaling approach and the Strichartz estimates, instead of the It\^o formula for $|Y(t)|_2^2$,
to control the singular coefficient $h_2(X^u)$
and to obtain pathwise estimates of solutions on small intervals, which then  by iteration yield the global pathwise estimates
\eqref{esti*}, \eqref{psi-Psi*} below.  To this aim, we consider in this case  the following basic hypothesis.

\begin{enumerate}
\item[(H2)] $2\leq \a < 1+ \frac{4}{d}$, $1\leq d\leq 3$, and $e_k$ are constants, $1\leq k\leq N$.
\end{enumerate}

(In the case where $e_k$ are not constant, which is ruled out here, there arise some delicate problems
related to the nonintegrability of $(B_j)^{c(B_j)^{2/\theta}}$, where $B_k:= \sup_{t\in [0,T]} |\beta_k(t)|$,
$\theta\in (0,1)$ and $c>0$, $1\leq k\leq N$.)

It is easily seen that $(H2)$ implies $(H0)$ and $(H1)$ and also that $(H2)$ is fulfilled in
some important physical models, for instance the Gross-Pitaevskii model when $d=1,2$ (\cite{11a}).
As a matter of fact, under Hypothesis $(H2)$, one has  not only \eqref{e2.6}
with equality, but also that the optimal pair $(X,u)$ satisfies
the stochastic maximum principle. The main result is formulated below.

\begin{theorem}\label{t2.6}
Consider  $\Phi$   with $\g_3=0$. Assume Hypothesis $(H2)$, and $\bbx_T\in L^{2+\nu}(\Omega;H^1) $,
$ \bbx_1 \in L^{2+ \nu}(\Omega; L^2(0,T; H^1)) $  for some small $\nu \in (0,1)$.

Then, for each $x\in H^1$, $0<T<\9$, there exists a relaxed solution
$(\Omega^*, \mathscr{F}^*, \bbp^*, (\mathscr{F}^*_t)_{t\geq 0}, \beta^*, u^*, X^*)$
in the sense of Definition \ref{def-Phi} to Problem (P), such that
\begin{equation}\label{e2.10}
 \Phi^*(u^*)=\inf\{\Phi(u);\ u\in\calu_{ad},\ X^u\ satisfies\ \eqref{equa-x} \}.
\end{equation}
Moreover, we have (the stochastic maximum principle)
\begin{equation}\label{e5.1}
u^* (t)=P_U\(\frac 1{\g_2}\,{\rm Im}\int_{\rr^d}V(\xi) X^*(t,\xi) \ol{Y^*}(t,\xi)d\xi\),\ \ff t\in[0,T],\ \bbp^*-a.s.
\end{equation}
where $P_U$ is the projection on $U$, and $(Y^*,Z^*)$ is the solution to the dual backward  stochastic  equation \eqref{back-equa}
with
$\bbx_T, \bbx_1, \beta, u, X^u$ replaced by $\bbx^*_T, \bbx^*_1, \beta^*, u^*, X^*$ respectively.
\end{theorem}

In the deterministic case (i.e. $\mu_k=0$, $1\leq k\leq N$), for the initial datum $x\in H^1$,
the optimal control indeed exists
for the exponent $\a\geq 2$ and in the
energy-subcritical case $(H1)$, which is also new in the literature.

\begin{theorem} \label{Thm-deter}
In the deterministic case (i.e. $\mu_k=0$, $1\leq k\leq N$),
consider $\Phi$ in \eqref{def-Phi} with $\g_3=0$ and the exponent $\a\geq 2$ in the range specified in Hypothesis $(H1)$.
Assume that $\bbx_T\in H^1 $ and
$ \bbx_1 \in  L^2(0,T; H^1) $.

Then, for each $x\in H^1$, $0<T<\9$, there exists an optimal control $u$ to Problem (P) such that
\begin{equation}\label{e2.10-deter}
\Phi(u)=\inf\{\Phi(v);\ v\in\calu_{ad},\ X^v\ satisfies\ \eqref{e2.1}\}.
\end{equation}
Moreover,
\begin{equation}\label{e5.1-deter}
u(t)=P_U\(\frac 1{\g_2}\,{\rm Im}\int_{\rr^d}V(\xi) X(t,\xi) \ol{Y}(t,\xi)d\xi\),\ \ff t\in[0,T],
\end{equation}
where $P_U$ is the projection on $U$, and $Y$ is the solution to the backward   equation \eqref{back-equa} with $Z=0$.
\end{theorem}

\begin{remark}
Optimal bilinear control is also studied in \cite{11} and \cite{11a} for  linear and nonlinear deterministic
Schr\"odinger equations respectively. In both papers, some compactness conditions of initial data or  controls are needed
for the existence of the optimal control.
More precisely, in \cite{11} the initial data belong to a compact subspace of $L^2$, while in \cite{11a} the minimizing
controls are bounded in $H^1[0,T]$, hence compact in $L^2[0,T]$. In contrast to this, in Theorem \ref{Thm-deter},
the existence of the optimal control is obtained without these conditions, and the proof is quite different and applies as well to the stochastic
case.
Moreover, unlike in \cite{11a}, less regularity of the initial data is required in  Theorem  \ref{Thm-deter}
for the maximum principle \eqref{e5.1-deter}.
\end{remark}

The proof of Theorem \ref{Thm-deter} follows the lines of that of Theorem \ref{t2.6}, the main part of which are
the analysis of equation of variation \eqref{forward-equa} as well as of the backward stochastic equation \eqref{back-equa},
and the tightness of controls.

The key idea to obtain the tightness of controls in
this irregular case is
to employ the Ekeland principle, as well as the directional derivative of $\Phi$, to obtain
the representation formula of the minimizing controls (see \eqref{e4.13} below).
Then, by virtue of the integrability
of the forward and  backward solutions to  \eqref{equa-x} and \eqref{back-equa} respectively, one is able to obtain
the tightness of controls
in the space $L^1(0,T; \bbr^d)$, which consequently yields
equality in \eqref{e2.6} by analogous arguments as in the proof of Theorem \ref{t2.5}.

As mentioned above, the proof of integrability of the stochastic backward solution relies on duality analysis, which
is also of independent interest.

The remaining part of this paper is organized as follows. Section \ref{PROOF-THM1} includes the proof of Theorem \ref{t2.5}.
Section \ref{GD-PHI} and Section \ref{PROOF-THM2} are mainly devoted to the proof of Theorem \ref{t2.6}.
Section \ref{GD-PHI} is concerned with the directional derivative of $\Phi$, which requires the analysis of the equation
of variation \eqref{forward-equa} and of the backward stochastic equation \eqref{back-equa}. Section \ref{PROOF-THM2}
mainly contains the proof for the tightness of controls. The proof of Theorem \ref{Thm-deter} is also included there.
For simplicity of the exposition, some auxiliary lemmas and technical proofs are postponed to the Appendix, i.e. Section \ref{APPDIX}.

\paragraph{Notations.} For $1\le p\le\9$, we denote by $L^p(\rr^d)=L^p$ the space of all Lebesgue $p$-integrable (complex-valued) functions on the real Euclidean space $\rr^d$.
The norm of $L^p$ is denoted by $|\cdot|_{L^p}$, and $p'\in[1,\9]$ denotes the unique number such that $\frac1p+\frac1{p'}=1$. In particular, the Hilbert space $L^2(\rr^d)$ is endowed with the scalar product
$$\<y,z\>_2=\int_{\rr^d} y(\xi)\bar z(\xi)d\xi;\ \ y,z\in L^2(\rr^d),$$where $\bar z$ is the complex conjugate of $z\in\mathbb{C}$. We also use $|\cdot|_2=|\cdot|_{L^2}$.
$W^{1,p}=W^{1,p}(\rr^d)$ is the classical Sobolev space $\{v\in L^p;\ \nabla v\in L^p\}$ with the norm $\|v\|_{W^{1,p}} = |v|_2 + |\na v|_2$, $H^1=W^{1,2}$ and $H^{-1}$ is the dual space of $H^1$.

By $L^q(0,T;L^p)$ we denote the space of all integrable $L^p$-valued functions $u:(0,T)\to L^p$ with the norm
$$\|u\|_{L^q(0,T;L^p)}=\left(\int^T_0\left(\int_{\rr^d}|u(t,\xi)|^pd\xi\right)^{\frac qp}dt\right)^{\frac1q}.$$
By $C([0,T];L^p)$ we denote the standard space of all $L^p$-valued continuous functions on $[0,T]$ with the sup norm in $t$.
$L^q(0,T; W^{1,p})$ and $C([0,T]; H^1)$  are defined similarly.
$\mathcal{D}(0,T; \bbr^m)$ is the set of all $\bbr^m$-valued  smooth and compactly supported functions,
and $\mathcal{D}'(0,T; \bbr^m)$ is its dual space.

We denote by $|\cdot|_m$ the Euclidean norm in $\rr^m$ and by $u\cdot v$ the scalar product of vectors $u,v\in\bbr^m$.
We shall use the standard notations to represent spaces of infinite dimensional stochastic processes (see, e.g., \cite{6, LR15}).

Throughout this paper, we use $C$ for various constants that may change from line to line.

\section{Proof of Theorem \ref{t2.5}. } \label{PROOF-THM1}

We set
\begin{align*}
   I:= \inf \{\Phi(u); u\in \mathcal{U}_{ad}, \ X^u\ satisfies\ \eqref{equa-x}\}> 0
\end{align*}
and consider a sequence $\{u_n\} \subset \mathcal{U}_{ad}$ such that
\begin{align} \label{bbd-Psi}
      I \leq \Phi(u_n) \leq I + n^{-1}, \ \ \forall n\in \mathbb{N}.
\end{align}
Since $\g_3 >0$, this yields
\begin{align} \label{bdd-un}
   \sup\limits_{n\geq 1} \bbe \int_0^T (|u_n(t)|_m^2 + |u'_n(t)|_m^2) dt  <\9.
\end{align}

\begin{lemma} \label{Lem-tight}
Let $\bbp \circ u_n^{-1}$ be the probability measures induced by $u_n$, $n\geq 1$. Then,
$\{\bbp \circ u_n^{-1}\}$ is tight in $C([0,T]; \bbr^m)$.
\end{lemma}

{\it \bf Proof.}  By the Arzel\`{a} theorem, it suffices to show that
\begin{align} \label{AA-bdd}
   \lim\limits_{R\to \9} \sup\limits_{n\geq 1}
   \bbp \circ u_n^{-1}\bigg\{v\in C([0,T];\bbr^m): \sup\limits_{t\in[0,T]} |v(t)|_m > R \bigg\} = 0,
\end{align}
and for any $\ve >0$,
\begin{align} \label{AA-equicon}
    \lim\limits_{\delta\to 0} \sup\limits_{n\geq 1}
    \bbp \circ u_n^{-1} \bigg\{v\in C([0,T]; \bbr^m): \sup\limits_{|t-s|\leq \delta}  |v(t)-v(s)|_m >\ve \bigg\} =0.
\end{align}

In fact, \eqref{AA-bdd} follows immediately form the uniform boundedness of $\{u_n\}$,
while \eqref{AA-equicon} follows by \eqref{bdd-un} and
\begin{align}
  &\sup\limits_{n\geq 1} \bbp \circ u_n^{-1} \bigg\{v\in C([0,T]; \bbr^m): \sup\limits_{|t-s|\leq \delta}  |v(t)-v(s)|_m >\ve \bigg\} \nonumber \\
  \leq& \frac 1 \ve \sup\limits_{n\geq 1} \bbe \sup\limits_{|t-s|\leq\delta} |u_n(t)-u_n(s)|_m \nonumber \\
  \leq& \frac {\delta^{\frac 12}} {\ve}  \sup\limits_{n\geq 1} \bbe  \|u'_n\|_{L^2(0,T;\bbr^m)} \to 0, \ \ as\ \delta\to 0.
\end{align}
 \hfill $\square $

Now, consider the sequence $\calx_n:=(\bbx_T, \bbx_1, \beta, u_n)$ with $\beta=(\beta_1,\cdots,\beta_N)$ in the space
$\mathcal{Y}:=   L^2(  \bbr^d) \times L^2( (0,T)\times \bbr^d) \times C([0,T];\bbr^N) \times C([0,T]; \bbr^m)$.
Lemma \ref{Lem-tight} implies that the induced probability measures of $\calx_n$, $n\in\mathbb{N}$, are tight in the space $\mathcal{Y}$. Then,
by Prohorov's theorem, they are weakly compact and so, by the Skorohod representation theorem, there exist a probability space
$(\Omega^*, \mathcal{F}^*, \bbp^*)$ and  $\calx^*_n := ((\bbx^*_{T})_n, \bbx^*_{1,n},\beta^*_n, u^*_n)$,
$\calx^* :=(\bbx^*_T,\bbx^*_{1}, \beta^*, u^*)$ in $\mathcal{Y}$, $n\in \mathbb{N}$, such that the joint distribution of $\calx_n^*$ and $\calx_n$
coincide,  and $\bbp^*$-a.s., as $n\to \9$,
\begin{align}
      &\beta^*_{n} \to \beta^*\ \ in\ C([0,T]; \bbr^N), \label{conv-bn*}\\
      &(\bbx^*_{T})_n \to \bbx^*_T,\ \ in\ L^2(\bbr^d),\ \   \bbx^*_{1,n} \to \bbx^*_{1},\ \ in\ L^2( (0,T)\times \bbr^d), \label{conv-XT*-X1*}
\end{align}
and
\begin{align}
      &u^*_n \to u^*\ \ in\ C([0,T]; \bbr^m). \label{conv-un*}
\end{align}
Note that, $((\bbx_T^*)_n, \bbx^*_{1,n}, \beta^*_n)$ has the same distribution as $(\bbx_T, \bbx_1, \beta)$,
and so does the limit $(\bbx^*_T, \bbx^*_1, \beta^*)$.

For each $n\geq 1$, define $\calf^*_{t,n}:=\sigma(\calx^*_n(s), s\leq t)$. Then,
$\bbx^*_n(T) \in L^2(\Omega, \mathcal{F}^*_T, \bbp^*; L^2)$, $\bbx^*_{1,n} \in L^2_{ad^*}(0,T; L^2(\Omega; L^2))$,
$u^*_n \in \calu_{ad^*}$, and
$(\beta^*_n(t), \calf^*_{t,n})$, $t\in [0,T]$, is a Wiener process.
It follows from Proposition $2.2$ that, under the hypothesis $(H0)$, for each $(\beta^*_n, u^*_n)$ there exists
a unique strong $L^2$-solution $X^*_n$ to \eqref{equa-x}. Hence, $(\Omega^*, \mathcal{F}^*, \bbp^*, (\mathcal{F}^*_t)_{t\geq 0}, \beta_n^*, u_n^*, X_n^*)$
is an admissible system.

Moreover, since the solution to \eqref{equa-x} is a measurable map of Brownian motions and controls,
we also have that the distribution of $((\bbx^*_T)_n, \bbx^*_{1,n}, \beta_n^*, u_n^*, X_n^*)$ are the same to that of
$(\bbx_T, \bbx_{1}, \beta, u_n, X_n)$, where $X_n$ is the solution to \eqref{equa-x} corresponding to $(\beta, u_n)$.
In particular, $\Phi^*(u^*_n) = \Phi(u_n)$.

Similarly, set $\calf^*_t:= \sigma(\calx^*(s), s\leq t)$ and let $X^*$ be the unique strong $L^2$-solution to \eqref{equa-x} corresponding to
$(\beta^*, u^*)$. Then, $(\Omega^*, \mathcal{F}^*, \bbp^*, (\mathcal{F}^*_t)_{t\geq 0}, \beta^*, u^*, X^*)$
is an admissible system.

Below, we consider the derivatives of $u^*_n$ and $u^*$.
For each $n\geq 1$,  define $(u_n^*)' \in \mathcal{D}'(0,T; \bbr^m)$  in the distribution sense, i.e.,
$((u_n^*)' , v)= - (u_n^*, v'), \forall v\in \mathcal{D}(0,T; \bbr^m)$, where
$(\ ,\ )$ denotes the pairing between $\mathcal{D}'(0,T; \bbr^m)$ and $\mathcal{D}(0,T; \bbr^m)$.
We claim that $(u_n^*)'$ has the same distribution as $u_n'$ in $\mathcal{D}'(0,T; \bbr^m)$,
and $\bbe^*\|(u^*_n)'\|^2_{L^2(0,T; \bbr^m)} = \bbe\|u_n'\|^2_{L^2(0,T; \bbr^m)}$.
It follows from  \eqref{bdd-un} that there exists $v^* \in L^2(\Omega^*; L^2(0,T;\bbr^m))$, such that
\begin{align} \label{conv-v*}
     (u_n^*)' \to  v^*,\ \ weakly\ in \ L^2(\Omega^*; L^2(0,T; \bbr^m)),\ n\to \9.
\end{align}
Indeed, for any $l\geq 1$, $v_j \in \mathcal{D}(0,T; \bbr^m)$,
$1\leq j\leq l$, and $c \in \bbr^l$,
\begin{align*}
  &\bbe^* exp(i\sum\limits_{j=1}^l c_j ((u^*_n)', v_j))
  =  \bbe^* exp(-i(u^*_n, \sum\limits_{j=1}^l c_j  v'_j)) \\
  =&  \bbe exp(-i(u_n, \sum\limits_{j=1}^l c_j  v'_j))
  = \bbe exp(i\sum\limits_{j=1}^l c_j (u'_n, v_j)),
\end{align*}
which implies that
the distributions of $(u^*_n)'$ and $u'_n$ coincide.
Moreover, if $\mathscr{D}:= \{v_n\}$ is a dense subset in $\mathcal{D}(0,T; \bbr^m)$ (hence also dense in $L^2(0,T; \bbr^m)$),
we have
\begin{align*}
     \bbe^* \sup\limits_{n\geq 1} \frac{|((u^*_n)',v_n)|^2}{ \|v_n\|^2_{L^2(0,T; \bbr^m)}}
    = \bbe \sup\limits_{n\geq 1} \frac{|(u'_n,v_n)|^2}{ \|v_n\|^2_{L^2(0,T; \bbr^m)}}
    = \bbe \|u'_n\|^2_{L^2(0,T; \bbr^m)} < \9,
\end{align*}
which implies that $\bbe^*\|(u^*_n)'\|^2_{L^2(0,T; \bbr^m)} = \bbe\|u_n'\|^2_{L^2(0,T; \bbr^m)}$, as claimed.

Similarly, define $(u^*)' \in \mathcal{D}'(0,T; \bbr^m)$ in the distribution sense.
We have
\begin{align} \label{v*-u'*}
      (u^*)' = v^*,\ \ in\ L^2(0,T; \bbr^m),\ \bbp^*-a.s.
\end{align}
Indeed, let $\mathcal{E}$ be a countable dense set in  $L^\9(\Omega^*)$.
For any $v\in \mathscr{D}$ and $\psi \in \mathcal{E}$, by \eqref{conv-un*}, the dominated convergence theorem and the
weak convergence \eqref{conv-v*}, it follows that
\begin{align*}
  &\bbe^* ((u^*)',v) \psi
  = - \bbe^* \psi (\int_0^T u^*(t) \cdot v'(t)dt)
  = - \lim\limits_{n\to \9}  \bbe^* \psi(\int_0^T u_n^*(t) \cdot v'(t)dt) \\
  =&  \lim\limits_{n\to \9} \bbe^* \int_0^T (u_n^*)'(t) \cdot v(t) \psi dt
  = \bbe^* \int_0^T v^*(t) \cdot v(t)\psi dt
  = \bbe^* (v^*,v) \psi.
\end{align*}
Hence,  $ ((u^*)',v) = (v^*,v)$, $\bbp^*$-a.s., $v \in \mathscr{D}$.  Since
$\mathscr{D}$ is countable and dense in $L^2(0,T; \bbr^m)$,  \eqref{v*-u'*} follows. \\

Next, we show that the solutions to \eqref{equa-x} depend pathwisely continuous with respect to controllers,
by using the  rescaling approach developed in \cite{1,2}.

\begin{lemma} \label{Lem-Xn*}
Let $X^*_n$ (resp. $X^*$) be the solution to \eqref{equa-x} corresponding to $(\beta^*_n, u^*_n)$ (resp. $(\beta^*, u^*)$) as above, $n\geq 1$.
Assume the conditions in Theorem $2.5$ to hold. Then, for each $x\in L^2$, $0<T<\9$ and Strichartz pair $(p,q)$,    we have $\bbp^*$-a.s., as
$n\to \9$,
\begin{align} \label{conv-Xn*}
   \|X^*_n - X^*\|_{C([0,T];L^2)} +  \|X^*_n - X^*\|_{L^q(0,T; L^p)} \to 0.
\end{align}
\end{lemma}

{\it \bf Proof. } Set $W^*_n(t,\xi) = \sum_{j=1}^N \mu_j e_j(\xi) \beta^*_{j,n}(t)$,
$W^*(t,\xi) = \sum_{j=1}^N \mu_j e_j(\xi) \beta^*_{j}(t)$, $t\geq 0$, $\xi\in \bbr^d$.
We may assume $T \geq 1$ without loss of generality.

Note that,
in the conservative case,
\begin{align} \label{bdd-Xn*}
     |X^*_n(t)|_2 = |x|_2 < \9, \ \ t\in[0,T],\ n\geq 1,\  \bbp^*-a.s.
\end{align}

Using the rescaling transformation,
\begin{align} \label{rescal}
     y^*_n = e^{-W^*_n} X^*_n,
\end{align}
we deduce from \eqref{equa-x} with $X$, $u$, $\beta_j$ replaced by $X^*_n$, $u^*_n$ and $\beta_{j,n}^*$ respectively that
\begin{align} \label{equa-yn*}
  dy^*_n
  =& A^*_n(t)y^*_n dt - \lbb i |y^*_n|^{\a-1}y^*_n dt,
      +  f(u_n^*) y^*_n dt  \\
  y^*_n(0)=&x, \nonumber
\end{align}
where $A^*_n(t) = -i (\Delta + b^*_n(t)\cdot \na + c^*_n(t))$,
$b^*_n(t)= 2 \na W^*_n(t)$,
$c^*_n(t)= \sum_{j=1}^d (\partial_j W^*_n(t))^2 + \Delta W^*_n(t)$,
$\wt{\mu}=2^{-1}\sum_{j=1}^N \mu_j^2 e^2_j$, and
$f(u_n^*) = - i (  V_0 +  u^*_n \cdot V)$.

It suffices to prove that for each Strichartz pair $(p,q)$, $\bbp^*$-a.s.,
\begin{align} \label{conv-y-Xn*}
   \|y^*_n - y^*\|_{L^\9(0,T;L^2)} +  \|y^*_n - y^*\|_{L^q(0,T; L^p)} \to 0,\ \ as\ n\to\9.
\end{align}
Now, we will prove \eqref{conv-y-Xn*} for the Strichartz pair  $(p,q)=(\a+1,\frac{4(\a+1)}{d(\a-1)})$.
The general case will follow immediately from the Strichartz estimates (see, e.g., \cite[Lemma 4.1]{1}).

To this end, we prove first  that
\begin{align} \label{bdd-y-un*}
   \sup\limits_{n\geq 1} \|y^*_n\|_{L^q(0,T; L^p)}  < \9,\ \ \bbp^*-a.s.
\end{align}

Applying the Strichartz estimates to \eqref{equa-yn*} yields
\begin{align*}
     \|y^*_n\|_{L^q(0,t; L^p)}
    \leq  C_T \bigg[|x|_2
          + \|\lbb i |y^*_n|^{\a-1}y^*_n\|_{L^{q'}(0,t; L^{p'})}
        + \|f(u^*_n) y^*_n\|_{L^1(0,t;L^2)} \bigg].
\end{align*}
(Note that, the Strichartz coefficient $C_T$ is independent of $n$, since by \eqref{conv-bn*} $\sup_{n\geq 1}\|W^*_n\|_{C([0,T];L^\9)} <\9$, $\bbp^*$-a.s.)

By H\"older's inequality,
\begin{align} \label{esti-lp}
  \|\lbb i  |y^*_n|^{\a-1}y^*_n\|_{L^{q'}(0,t; L^{p'})}
  \leq t^\theta \|y^*_n\|^\a_{L^q(0,T;L^p)},
\end{align}
where $\theta= 1- \frac{d(\a-1)}{4} >0$, and by the conservation \eqref{bdd-Xn*},
\begin{align*}
  \|f(u^*_n)y^*_n\|_{L^1(0,t;L^2)}
  \leq T( |V_0|_{L^\9} + D_U\|V\|_{L^\9(\bbr^d;\bbr^m)})|x|_2.
\end{align*}
Thus,
\begin{align}\label{esti-yn*}
     \|y^*_n\|_{L^q(0,t; L^p)}  \leq C_T (D(T)|x|_2 + t^\theta \|y^*_n\|^\a_{L^q(0,T; L^p)}),
\end{align}
where $ D(T):= 1+T ( | V_0|_{\9}
          + D_U \|V\|_{L^\9(0,T; \bbr^m)} ). $

Choose $t \in [0,T]$ such that
$ C_T D(T) (|x|_2+1) = (1-\frac 1 \a) (\a C_T t^\theta)^{-\frac{1}{\a-1}}$,
i.e.,
\begin{align*}
    t = \a^{-\frac{\a}{\theta}} (\a-1)^{\frac{\a-1}{\theta}} (|x|_2+1)^{-\frac{\a-1}{\theta}} C_T^{-\frac{\a}{\theta}} D(T)^{-\frac{\a-1}{\theta}} (\leq T).
\end{align*}
Then, by Lemma \ref{Lem-Bound}, we get
\begin{align*}
    \|y^*_n\|_{L^q(0,t; L^{p})} \leq \frac{\a}{\a-1} C_T D(T) |x|_2 .
\end{align*}
Iterating similar estimates on $[jt,(j+1)t\wedge T]$, $1\leq j\leq [\frac{T}{t}]$,
we obtain
\begin{align} \label{esti-yn*-lpq}
     \|y^*_n\|_{L^q(0,T; L^p)}
     \leq& \([\frac{T}{t}]+1\)^{\frac 1 q} \frac{\a}{\a-1} C_T D(T) |x|_2 ,
\end{align}
which yields \eqref{bdd-y-un*}. \\

It remains to prove \eqref{conv-y-Xn*}. Applying the Strichartz estimates to the equations of $y^*_n$ and $y^*$, we have for any $t\in(0,T)$,
\begin{align*}
       & \|y^*_n - y^*\|_{L^\9(0,t;L^2)} +   \|y^*_n - y^*\|_{L^q(0,t; L^p)}  \\
       \leq& C_T \bigg[\|f(u_n^*) y^*_n - f(u^*) y^*\|_{L^1(0,t;L^2)}
            +   \|  |y^*_n|^{\a-1}y^*_n -    |y^*|^{\a-1}y^*  \|_{L^{q'}(0,t; L^{p'})} \bigg].
\end{align*}
Proceeding as in \cite[(4.12)]{1},  we have
\begin{align*}
    &\|  |y^*_n|^{\a-1}y^*_n -   |y^*|^{\a-1}y^* \|_{L^{q'}(0,t; L^{p'})} \\
    \leq& \a t^\theta (\|y^*_n\|^{\a-1}_{L^q(0,T;L^p) } + \|y^*\|^{\a-1}_{L^q(0,T;L^p)}) \| y^*_n - y^*\|_{L^q(0,t;L^p) }.
\end{align*}
Moreover, for $t\leq 1$,
\begin{align*}
     & \|f(u^*_n)y^*_n- f(u^*)y^*\|_{L^1(0,t;L^2)}
     \leq  t^\frac 12  D \( \|u^*_{n} - u^* \|_{L^2(0,T;\bbr^m)}
          +  \|y^*_n-y^*\|_{L^\9(0,t; L^2)} \),
\end{align*}
where  $D = 2(|V_0|_{L^\9} + \|V\|_{L^\9(\bbr^d; \bbr^m)} (|x|_2+D_U))$.

Hence,
\begin{align} \label{diff-yn*}
     & \|y^*_n - y^*\|_{L^\9(0,t;L^2)} +   \|y^*_n - y^*\|_{L^q(0,t; L^p)} \nonumber  \\
     \leq& C_T\bigg[ \a t^\theta  (\|y^*_n\|^{\a-1}_{L^q(0,T;L^{p})} + \|y^*\|^{\a-1}_{L^q(0,T;L^{p})} ) \|y^*_n - y^*\|_{L^{q}(0,t; L^p)} \nonumber \\
     & \qquad  + t^{\frac 12}  D \|u^*_{n} - u^*\|_{L^2 (0,T;\bbr^m)}
               + t^{\frac 12} D \|y^*_n - y^*\| _{L^\9(0,t; L^2)} \bigg] \nonumber \\
     \leq& (t^\theta + t^{\frac 12})C_T \wt{D}(T) (\|y^*_n- y^*\|_{L^\9(0,t;L^{2})}+\|y^*_n- y^*\|_{L^q(0,t;L^{p})} +\|u^*_{n} - u^*\|_{L^2 (0,T;\bbr^m)}),
\end{align}
where $ \wt{D}(T):=  \a (\sup _{n\geq 1}\|y^*_n\|^{\a-1}_{L^q(0,T;L^{p})} + \|y^*\|^{\a-1}_{L^q(0,T;L^{p})} ) +D <\9$,
$ \bbp^*$-a.s. Choosing $t$ small enough and independent of $n$, such that $t^\theta + t^{\frac{1}{2}} \leq (2\wt{D}(T)C_T)^{-1}$, we
get that $\bbp^*$-a.s. as $n\to \9$,
\begin{align} \label{esti-yn*-diff}
       \|y^*_n - y^*\|_{L^\9(0,t;L^2)} +   \|y^*_n - y^*\|_{L^q(0,t; L^p)}
     \leq   2  \|u^*_{n} - u^*\|_{L^2 (0,T;\bbr^m)} \to 0.
\end{align}
Since $t$ is independent of the initial data, iterating this procedure finite times we obtain  \eqref{conv-y-Xn*}, thereby completing the proof.
\hfill $\square $

\begin{lemma} \label{Lem-conv}
Let $X^*_n$, $(\bbx^*_T)_n$, $\bbx^*_{1,n}$, $X^*$, $\bbx^*_T$ and $\bbx^*_1$ be as above, $n\geq 1$.
We have $\bbp^*$-a.s., as $n\to \9$,
\begin{align} \label{App-conv-Xtau}
      \bbe^* Re\<X^*_n(T), (\bbx^*_T)_n\>_2
     \to \bbe^* Re \<X^*(T), \bbx^*_T\>_2,
\end{align}
and
\begin{align} \label{App-conv-intX}
     \bbe^* \int_0^{T} Re\<X^*_n(t), \bbx^*_{1,n}(t)\>_2 dt
     \to \bbe^* \int_0^{T} Re\<X^*(t), \bbx^*_{1}(t)\>_2 dt.
\end{align}
\end{lemma}

{\it \bf Proof.} By \eqref{conv-Xn*} and \eqref{conv-XT*-X1*}, $\bbp^*$-a.s., as $n\to \9$,
\begin{align} \label{App-conv-X.1}
      Re \<X^*_n(T), (\bbx^*_T)_n\>_2
     \to Re \<X^*(T), \bbx^*_T\>_2,
\end{align}
\begin{align}   \label{App-conv-X.2}
      \int_0^{T} Re \<X^*_n(t),\bbx^*_{1,n}(t)\>_2 dt
     \to \int_0^{T} Re \<X^*(t), \bbx^*_{1}(t)\>_2 dt.
\end{align}

Then, for any $\ve \in (0,1)$ fixed, by the Young inequality $ab\leq \frac{1-\ve}{2}a^{\frac{2}{1-\ve}} + \frac{1+\ve}{2} b^{\frac{2}{1+\ve}}$,
we get
\begin{align*}
    &\sup\limits_{n\geq 1} \bbe^* |\<X^*_n(T), (\bbx^*_T)_n\>_2|^{1+\ve}
    \leq  \sup\limits_{n\geq1} \bbe^* |X^*_n(T)|_2^{1+\ve} |(\bbx^*_T)_n|_2^{1+\ve} \\
    \leq&  \frac{1-\ve}{2} \sup\limits_{n\geq 1}\bbe^* |X^*_n(T)|_2^{\frac{2(1+\ve)}{1-\ve}}
           + \frac{1+\ve}{2} \sup\limits_{n\geq 1} \bbe^* |(\bbx^*_T)_n|_2^2 \\
    =&  \frac{1-\ve}{2}   |x|_2^{\frac{2(1+\ve)}{1-\ve}}
           + \frac{1+\ve}{2}\bbe  | \bbx_T |_2^2,
\end{align*}
which implies the uniform integrability of $\{Re \<X^*_n(T), (\bbx^*_T)_n\>_2\}_{n\geq 1}$, thereby
yielding \eqref{App-conv-Xtau} by \eqref{App-conv-X.1}.

Similarly, for $\ve \in (0,1)$ fixed, we have
\begin{align*}
    &\sup\limits_{n\geq 1}\bbe^* \bigg| \int_0^{T}  Re \<X^*_n(t), \bbx^*_{1,n}(t)\>_2 dt \bigg|^{1+\ve} \\
   \leq& \frac{1-\ve}{2} \sup\limits_{n\geq 1} \bbe^*   \|X^*_n\|_{L^2(0,T; L^2)}^{\frac{2(1+\ve)}{1-\ve}}
          + \frac{1+\ve}{2} \sup\limits_{n\geq 1} \bbe^*  \|\bbx^*_{1,n}\|_{L^2(0,T; L^2)}^2    \\
    = &  \frac{1-\ve}{2} T^{\frac{1+\ve}{1-\ve}} |x|_2^{\frac{2(1+\ve)}{1-\ve}}
          +\frac{1+\ve}{2} \bbe \int_0^{T} |\bbx_1(t)|_2^2 dt <\9,
\end{align*}
which in view of \eqref{App-conv-X.2} implies \eqref{App-conv-intX}, as claimed. \hfill $\square$ \\

{\it \bf Proof of Theorem \ref{t2.5}.} By the conservation identity \eqref{bdd-Xn*} we have
\begin{align*}
        \Phi^*(u^*_n)
     =& (1+\g_1)|x|_2^2 + \bbe^*|(\bbx^*_T)_n|_2^2 + \g_1 \bbe^* \int_0^{T} |\bbx^*_{1,n}(t)|_2^2 dt\\
      & -2 \bbe^* Re \<X^*_n(T), (\bbx^*_T)_n\>_2  - 2 \g_1 \bbe^* \int_0^{T} Re \<X^*_n(t), \bbx^*_{1,n}(t)\>_2 dt \\
      &  + \g_2 \bbe^* \int_0^{T} |u^*_n(t)|_m^2 dt +  \g_3 \bbe^* \int_0^{T} |(u^*_n)'(t)|_m^2 dt,
\end{align*}

Note that, since the distributions of $(\bbx^*_T, \bbx^*_1)$ and  $((\bbx^*_T)_n, \bbx^*_{1,n})$ coincide for
$n\geq 1$, we have
\begin{align} \label{X-X1}
   \bbe^*| \bbx^*_T |_2^2 + \g_1 \bbe^* \int_0^{T} |\bbx^*_{1}(t)|_2^2 dt
  = \lim \limits_{n\to \9} \(\bbe^*|(\bbx^*_T)_n|_2^2 + \g_1 \bbe^* \int_0^{T} |\bbx^*_{1,n}(t)|_2^2 dt\).
\end{align}
Moreover, by \eqref{conv-un*} and the bounded dominated convergence theorem, it follows that
\begin{align*}
   \g_2 \bbe^* \int_0^{T} |u^*_n(t)|_m^2 dt \to \g_2 \bbe^* \int_0^{T} |u^* (t)|_m^2 dt,\ as\ n\to \9,
\end{align*}
and by \eqref{conv-v*} and \eqref{v*-u'*},
\begin{align*}
      \bbe^* \int_0^T |(u^*)'(t)|^2_m dt
    \leq& \liminf\limits_{n\to \9}   \bbe^* \int_0^T |(u^*_n)'(t)|^2_m dt.
\end{align*}
Thus,  taking into account Lemma \ref{Lem-conv}, we obtain
\begin{align*}
     \Phi^*(u^*)
     \leq  \liminf\limits_{n\to\9} \Phi^*(u^*_n)
    =  \liminf\limits_{n\to \9} \Phi(u_n) = I.
\end{align*}
which completes the proof.
\hfill $\square$

\begin{remark}
The proofs  above show also that, in the case $\g_3=0$, the objective functional $\Phi$ depends continuously on controls.
\end{remark}

\section{The directional derivative of function $\Phi$} \label{GD-PHI}
\setcounter{equation}{0}

This section is devoted to the calculation of the directional derivative of function $ \Phi $ on the convex set $\calu_{ad}$.
Namely, one has

\begin{proposition}\label{l4.1}
Assume that $\g_3=0$ and that the conditions of Theorem \ref{t2.6} to hold.
Then, for each $x\in L^2$ and all $u, v\in \calu_{ad}$, we have
\begin{equation}\label{e4.1}
\lim_{\vp\to0}\frac1\vp\,(\Phi(u+\vp\wt u)-\Phi(u))=\E\int^{T}_0\eta(u)(t)\cdot{\wt u}(t)dt,
\end{equation}
where $\wt{u} = v-u$, and
\begin{equation}\label{e4.2}
\eta(u)=2\(\g_2 u-{\rm Im}\int_{\rr^d}V(\xi) X^u(\xi) \ol{Y^u}(\xi)d\xi\).
\end{equation}
Here $(Y^u, Z^u)$ is the solution to the dual backward stochastic equation \eqref{back-equa}.
\end{proposition}

To prove Proposition \ref{l4.1}, we first study
the equation of variation associated with Problem (P), namely,

\begin{align} \label{forward-equa}
  &id \vf =   \Delta \vf dt  + \lbb   h_1(X^u) \vf dt  + \lbb  h_2(X^u) \ol{\vf} dt -i \mu \vf dt  \nonumber \\
  &\qquad \quad  + V_0 \vf dt + u\cdot V \vf dt + \wt{u}\cdot V X^u dt +i \vf d W(t), \nonumber \\
  &\vf(0)= 0,
\end{align}
where $\wt{u}= v-u$, $u,v\in \calu_{ad}$, $X^u$ is the solution to \eqref{equa-x}, and $h_j(X^u)$, $j=1,2$, are defined as in
\eqref{e4.4}.
The strong $H^1$-(and $L^2$-)solution to \eqref{forward-equa} can be defined similarly as in Definition \ref{d2.1}.

\begin{lemma} \label{WP-Forward}
$(i)$  Under Hypothesis $(H0)$, for   $u,v \in \calu_{ad}$, $\wt{u}:=v-u$,
there exists a unique strong $L^2$-solution $\vf^{u,\wt{u}}$ to \eqref{forward-equa} on $[0,T]$.

$(ii)$ Under Hypothesis $(H2)$,
for
any Strichartz pair $(p,q)$,
\begin{align} \label{bdd-forward-y}
     \sup\limits_{u,v\in \calu_{ad}} \(\|\vf^{u,\wt{u}}\|_{L^\9(\Omega; C([0,T]; L^2))} + \|\vf^{u,\wt{u}}\|_{L^\9(\Omega;L^q(0,T; L^p))} \) < \9.
\end{align}

Moreover, set $u_\ve:= u+\ve \wt{u}$ and let
$X^u$ and  $X^{u_\ve}$ be the corresponding solutions to \eqref{equa-x} with the initial datum $x\in H^1$.  Then,
\begin{align} \label{asy-X-y}
    \lim\limits_{\ve \to 0} \bbe \sup\limits_{t\in [0,T]} |\ve^{-1}(X^{u_\ve}(t) - X^u(t)) -\vf^{u,\wt{u}}(t)|^2_2 =0.
\end{align}
\end{lemma}

\begin{remark} \label{Rem-Diff-h2}

In comparison with $(ii)$, the weaker Hypothesis $(H0)$ is sufficient for the pathwise existence and uniqueness of the solution to \eqref{forward-equa},
thanks to the linear structure of \eqref{forward-equa}.
However, as mentioned in Section \ref{PROBLEM-RESULT},
Hypothesis $(H2)$ is needed in order that the estimate \eqref{bdd-forward-y} holds. The arguments presented below, particularly in the proof of the estimate \eqref{bdd-forward-y},
will also be used in the  analysis of the dual backward stochastic equation in the proof of Proposition \ref{WP-Backward} below.

\end{remark}

{\it \bf Proof of Lemma \ref{WP-Forward}.}
$(i)$ We set $z^{u,\wt{u}}:= e^{-W} \vf^{u,\wt{u}}$, $\wt{u}:=v-u$, and
for simplicity, we omit the dependence of $u,\wt{u}$ in $z^{u,\wt{u}}$ below.
It follows from \eqref{forward-equa} that
\begin{align} \label{equa-z}
    &dz = A(t) z dt - \lbb i h(X^u, z) dt + f(u) z dt  - i \wt{u} \cdot V  e^{-W}X^u dt, \ t\in(0,T), \\
    &z(0)=0, \nonumber
\end{align}
where $A(t)$ is similar as $A^*_n(t)$ in \eqref{equa-yn*}, i.e.,
$A(t) = -i (\Delta + b(t)\cdot \na + c(t))$,
$b(t)= 2 \na W(t)$,
$c(t)= \sum_{j=1}^d (\partial_j W(t))^2 + \Delta W(t)$,
$h(X^u,z):= h_1(X^u) z + h_2(X^u) e^{-2 i Im W} \ol{z}$,
and $f(u):= -i( V_0   +  u\cdot V )$.

It is equivalent to prove the existence and uniqueness of the solution to \eqref{equa-z} (see the proof of \cite[Lemma 6.1]{1}).

To this purpose, we reformulate \eqref{equa-z} in the mild form as
\begin{align} \label{equa-mild-z}
     z(t)  = \int_0^t U(t,s) \big[-\lbb i h(X^u, z)(s) & + f(u(s)) z(s)  - i\wt{u}(s)\cdot V e^{-W(s)}X^u(s) \big] ds,
\end{align}
where $0\leq t\leq T$, and $U(t,s)$, $0\leq s,t\leq T$, are the evolution operators corresponding to the operator $A(t)$
(see  \cite[Lemma 3..3]{1}). Choose the Strichartz pair $(p,q) = (\a+1, \frac{4(\a+1)}{d(\a-1)})$. Define the operator $F$ on $C([0,T]; L^2) \cap L^q(0,T;L^p)$ by
\begin{align*}
   F(\phi)(t) :=& \int_0^t U(t,s) \big[-\lbb i h(X^u,\phi)(s) +f(u(s)) \phi(s)  - i\wt{u}\cdot V e^{-W(s)}X^u(s) \big] ds,
\end{align*}
where $0\leq t\leq T$, $\phi \in C([0,T]; L^2) \cap L^q(0,T;L^p)$.
Set $\calz^{\tau_1}_{M_1} :=\{ \phi \in C([0,\tau_1]; L^2) \cap L^q(0,\tau_1;L^p): \|\phi\|_{C([0,\tau_1]; L^2)} + \|\phi\|_{L^q(0,\tau_1;L^p)} \leq M_1\} $,
where $\tau_1$ and $M_1$ are two random variables to be determined later.

Note that,  by H\"older's inequality,  for any $\phi_j \in C([0,T]; L^2) \cap L^q(0,T;L^p)$, $j=1,2$,
\begin{align*}
     \|h(X^u, \phi_1) - h(X^u, \phi_2)\|_{L^{q'}(0,t; L^{p'})}
    \leq   \a  t^\theta  \|X^u\|^{\a-1}_{L^q(0,t; L^{p})} \|\phi_1-\phi_2\|_{L^q(0,t; L^p)},
\end{align*}
where $\theta =  1- d(\a-1)/4 \in (0,1)$, and
\begin{align*}
    \|f(u)(\phi_1-\phi_2)\|_{L^{1}(0,t; L^{2})}
    \leq t \|f(u)\|_{L^\9(0,t; L^\9)}   \|\phi_1-\phi_2\|_{C([0,t]; L^2)}.
\end{align*}
Then, let $R_1(t) =\a  t^\theta   \|X^u\|^{\a-1}_{L^q(0,t; L^{p})} + t \|f(u)\|_{L^\9(0,t; L^\9)} $, $t\in [0,T]$.
By the Strichartz estimates and the above estimates we get for any $t\in [0,T]$,
\begin{align} \label{esti-diff}
    \|F(\phi_1)-F(\phi_2)\|_{C([0,t]; L^2) \cap L^q(0,t;L^p) }
   \leq  C_t R_1(t)   \|\phi_1-\phi_2\|_{C([0,t]; L^2)\cap L^q(0,t; L^p)}.
\end{align}
Similarly, for $\phi \in C([0,T]; L^2) \cap L^q(0,T;L^p)$, $t\in [0,T]$,
\begin{align*}
     \|F(\phi)\|_{C([0,t]; L^2) \cap L^q(0,t;L^p) }
     \leq C_t R_1(t)  \|\phi\|_{C([0,t]; L^2)\cap L^q(0,t; L^p)} + C_t \|\wt{u}\cdot V X^u\|_{L^1(0,t; L^2)}.
\end{align*}

Setting $\tau_1 = \inf\{t\in[0,T]:C_t R_1(t) \geq \frac 12\}\wedge T$, $M_1 = 2 C_{\tau_1} \|\wt{u}\cdot V X^u\|_{L^1(0,\tau_1; L^2)}$,
it follows by \eqref{esti-diff} that  $F$ is a contraction map in $\calz_{M_1}^{\tau_1}$, implying that there exists $\wt{z}_1 \in \calz_{M_1}^{\tau_1}$ such that $F(\wt{z}_1)=\wt{z}_1$.
Setting $z_1(\cdot):= \wt{z}_1(\cdot \wedge \tau_1)$  and using similar arguments as in \cite{1}, we deduce that $z_1$ is $(\mathcal{F}_t)$-adapted,
continuous in $L^2$, and solves
\eqref{equa-z} on $[0,\tau_1]$, and $\|z_1\|_{C([0,\tau_1]; L^2) \cap L^q(0,\tau_1;L^p) } \leq M_1 $. We also note that $\tau_1 \geq \sigma_*$,  where
\begin{align} \label{def-sig}
    \sigma_* := \inf\{t\in [0,T];  Z(t) \geq  \frac 12 \} \wedge T
\end{align}
with $ Z(t):= t^\theta \a  C_T \|X^u\|^{\a-1}_{L^q(0,T; L^{p})} + t C_T \|f(u)\|_{L^\9(0,T; L^\9)}$, $t\in [0,T].$

Suppose that at the $n^{th}$-step ($n\geq 1$) we have an increasing sequence of stopping times $\{\tau_j\}_{j=0}^n$ and an $L^2$-valued continuous $(\mathcal{F}_t)$-adapted process $z_n$, which satisfy that $\tau_0 = 0$, $\tau_j - \tau_{j-1} \geq \sigma_*$, $1\leq j\leq n$, $z_n$  solves \eqref{equa-z} on $[0,\tau_n]$, $z_n (\cdot )= z_n( \cdot \wedge \tau_n)$, and
$$\|z_n\|_{C([0,\tau_n]; L^2) \cap L^q(0,\tau_n;L^p) }
\leq  \sum\limits_{j=1}^n (2 C_{\tau_n})^{n+1-j}\ \|\wt{u}\cdot V X^u\|_{L^1(\tau_{j-1}, \tau_j; L^2)}.$$

Set $\calz^{\sigma_n}_{M_{n+1}} =\{\phi\in C([0,\sigma_n]; L^2) \cap L^q(0,\sigma_n;L^p): \|\phi\|_{C([0,\sigma_n]; L^2)} + \|\phi\|_{L^q(0,\sigma_n;L^p)} \leq M_{n+1}\} $,
where $\sigma_n$ and $M_{n+1}$ are random variables to be determined later.
Define the operator $F_n$ on $C([0,T]; L^2) \cap L^q(0,T;L^p)$ by
\begin{align*}
   F_n(\phi)(t) :=& U(\tau_n+t, \tau_n) z_n(\tau_n) + \int_0^t U(\tau_n+t,\tau_n+s) \big[-\lbb i h(X^u(\tau_n+s),\phi(s)) \\
                & \qquad \quad +f(u(\tau_n+s))\phi(s)  - i\wt{u}(\tau_n+s)\cdot V e^{-W(\tau_n+s)}X^u(\tau_n+s) \big] ds,
\end{align*}
where $0\leq t\leq T$, $\phi\in C([0,T]; L^2) \cap L^q(0,T;L^p)$.

Similarly,  for any $\phi_j \in \calz^{\sigma_n}_{M_{n+1}}$, $j=1,2$,
\begin{align*}
        &\|F_n(\phi_1)-F_n(\phi_2)\|_{C([0,\sigma_n]; L^2) \cap L^q(0,\sigma_n;L^p) } \\
   \leq&  C_{\tau_n+\sigma_n} R_{n+1}( \sigma_n)   \|\phi_1-\phi_2\|_{C([0,\sigma_n]; L^2)\cap L^q(0,\sigma_n; L^p)},
\end{align*}
where $R_{n+1}(t) =   \a  t^\theta\|X^u\|^{\a-1}_{L^q(\tau_n,\tau_n+t; L^{p})} + t \|f(u)\|_{L^\9(\tau_n,\tau_n+t; L^\9)} $, $t\in [0,T-\tau_n]$,
while for  $\phi \in \calz^{\sigma_n}_{M_{n+1}}$, we have
\begin{align*}
      &\|F_n(\phi)\|_{C([0,\sigma_n]; L^2) \cap L^q(0,\sigma_n;L^p) }  \\
     \leq& C_{\tau_n+\sigma_n} |z_n(\tau_n)|_2
           + C_{\tau_n+\sigma_n} \|\wt{u}\cdot V X^u\|_{L^1(\tau_n,\tau_n+\sigma_n; L^2)} \\
         &+  C_{\tau_n+\sigma_n} R_{n+1}( \sigma_n)   \|\phi \|_{C([0,\sigma_n]; L^2)\cap L^q(0,\sigma_n; L^p)} \\
     \leq& \frac 12 \sum\limits_{j=1}^{n+1} (2 C_{\tau_n+\sigma_n})^{n+2-j}  \|\wt{u}\cdot V X^u\|_{L^1(\tau_{j-1},\tau_j; L^2)} \\
         & +  C_{\tau_n+\sigma_n} R_{n+1}( \sigma_n)   \|\phi\|_{C([0,\sigma_n]; L^2)\cap L^q(0,\sigma_n; L^p)}.
\end{align*}
Then, let $\sigma_n(t) := \inf\{t\in[0,T -\tau_n]: C_{\tau_n+t} R_{n+1}( t) \geq \frac 12\} \wedge (T-\tau_n)$, $\tau_{n+1} := \tau_n + \sigma_n$,  and
$M_{n+1} :=  \sum_{j=1}^{n+1} (2 C_{\tau_{n+1}})^{n+2-j}\ \|\wt{u}\cdot V X^u\|_{L^1(\tau_{j-1}, \tau_j; L^2)} $.
It follows that $\tau_{n+1} - \tau_n = \sigma_n \geq \sigma_*$,  $F_n$ is a contraction map in $\calz^{\sigma_n}_{M_{n+1}}$, and so there exists
$\wt{z}_{n+1} \in \calz^{\sigma_n}_{M_{n+1}}$ satisfying $F_n(\wt{z}_{n+1}) = \wt{z}_{n+1}$.
As in \cite{1}, letting
\begin{align*}
    z_{n+1}(t) = \left\{
                \begin{array}{ll}
                  z_n(t), & \hbox{$t\in [0,\tau_n]$;} \\
                  \wt{z}_{n+1}((t-\tau_n)\wedge \sigma_n), & \hbox{$t\in (\tau_n, T]$,}
                \end{array}
              \right.
\end{align*}
it follows that $z_{n+1}$ is continuous $(\mathcal{F}_t)$-adapted, satisfies
\eqref{equa-z} on $[0,\tau_{n+1}]$, $z_n (\cdot )= z_n( \cdot \wedge \tau_{n+1})$, and
$$\|z_n\|_{C([0,\tau_{n+1}]; L^2) \cap L^q(0,\tau_{n+1};L^p) }
\leq  \sum\limits_{j=1}^{n+1} (2 C_{\tau_{n+1}})^{n+2-j}\ \|\wt{u}\cdot V X^u\|_{L^1(\tau_{j-1}, \tau_j; L^2)}.$$

Iterating this procedure, since $\sigma_n \geq \sigma_*$, we see that after at most $[T/\sigma_*]+1$
steps the stopping time $\tau_n$ reaches $T$. Hence, $\bbp$-a.s. there exists a global solution (denoted by $z$) on $[0,T]$ which satisfies
\begin{align} \label{esti*}
    \|z\|_{C([0,T]; L^2) \cap L^q(0,T;L^p) }
\leq  \sum\limits_{j=1}^{[T/\sigma_*]+1} (2 C_T)^{[T/\sigma_*]+2-j}\ \|\wt{u}\cdot V X^u\|_{L^1(\tau_{j-1}, \tau_j; L^2)}.
\end{align}

As regards the uniqueness, given any two solutions $\vf_j$,  we set $z_j = e^{-W} \vf_j$, $j=1,2$. Then, similarly to \eqref{esti-diff},
for any $s,t \in (0,T)$, $s+t \leq T$,  we have
\begin{align*}
  &\|z_1 - z_2\|_{C([s,s+t]; L^2) \cap L^q(s,s+t;L^p)}  \\
  \leq& C_T (\a t^\theta \|X^u\|^{\a-1}_{L^q(0,T;L^p)} + t \|f(u)\|_{L^\9(0,T; L^\9)})
        \|z_1 - z_2\|_{C([s,s+t]; L^2) \cap L^q(s,s+t;L^p)} ,
\end{align*}
which implies that $z_1=z_2$ on $[s,s+t]$, $\bbp$-a.s., for $t$ sufficiently small and independent of $s$,
thereby yielding the uniqueness by the arbitrariness of $s$. \\

$(ii)$ Under Hypothesis $(H2)$, the Strichartz coefficient $C_T$ is now identically a deterministic constant. Moreover,
\eqref{Cons-bdd-Xn-Lpq} and
\eqref{def-sig} imply that $\sigma_*$ has a deterministic lower bound, namely,
\begin{align*}
    \sigma_* \geq t_* := \inf \{t\in[0,T]: Z^*(t) \geq T\} \wedge T, \ \ \bbp-a.s.,
\end{align*}
where $Z^*(t) := \a C (t^\theta+t) \sup_{u\in \calu_{ad}} (\|X^u\|^{\a-1}_{L^\9(\Omega; L^q(0,T; L^p))} +  \|f(u)\|_{L^\9(\Omega;L^{\9}(0,T; L^\9))})$.
Thus, taking into account \eqref{esti*} and the uniform boundedness of $u,v\in \calu_{ad}$, we obtain \eqref{bdd-forward-y}.

Now, set $\wt{X}^u_\ve : = \ve^{-1}(X^{u_\ve} - X^u) - \vf$ and $\wt{y}^u_\ve:= e^{-W} \wt{X}^u_\ve $.
We need to prove that
\begin{align} \label{asy-X-eta}
    \lim\limits_{\ve \to 0} \bbe \|\wt{y}^u_\ve\|^2_{C([0,T]; L^2)} =0.
\end{align}

To this purpose,  note that
\begin{align*}
    \ve^{-1} (u_\ve \cdot V X^{u_\ve} - u \cdot V X^u) = \wt{u} \cdot V X^u + u_\ve \cdot V (\wt{X}^u_\ve + \vf),
\end{align*}
and
\begin{align*}
    &\ve^{-1} ( |X^{u_\ve}|^{\a-1}X^{u_\ve} - |X^u|^{\a-1}X^u) - (h_1(X^u) \vf + h_2(X^u)\ol{\vf})  \\
    =& \(\int_0^1 h_1(X_{u,r,\ve}) dr\) \wt{X}^u_\ve  + \(\int_0^1 h_2(X_{u,r,\ve}) dr\) \ol{\wt{X}^u_\ve} \\
     & + \vf \int_0^1 (h_1(X_{u,r,\ve}) dr - h_1 (X^u)) dr   +  \ol{\vf} \int_0^1 (h_2(X_{u,r,\ve}) dr - h_2(X^u)) dr.
\end{align*}
where $X_{u,r,\ve}= X^u + r(X^{u_\ve}-X^u)$, $r\in [0,1]$.
For simplicity, set
$R_j(\ve):= \int_0^1 (h_j(X_{u,r,\ve}) - h_j (X^u)) dr$, $j=1,2$,
and $R(\ve,\vf):=  -i(\lbb R_1(\ve)\vf + \lbb R_2(\ve)\ol{\vf} + \ve\wt{u}\cdot V\vf)$.

Then, by \eqref{equa-x} and \eqref{forward-equa}, $\wt{X}^u_\ve$ satisfies the equation
\begin{align} \label{equa-Xue}
    d \wt{X}^u_\ve  =&  -i \Delta \wt{X}^u_\ve dt
                      - \lbb i \int_0^1 h_1(X_{u,r,\ve})  dr \wt{X}^u_\ve dt
                       - \lbb i \int_0^1 h_2(X_{u,r,\ve}) dr \ol{\wt{X}^u_\ve}  dt
                         \nonumber  \\
                     &  -(\mu+iV_0 + iu_\ve \cdot V) \wt{X}^u_\ve dt + R (\ve, \vf)  dt
                        +  \wt{X}^u_\ve  d W(t).
\end{align}
This yields
\begin{align} \label{equa-yue}
      d \wt{y}^u_\ve
      = -i\Delta \wt{y}^u_\ve dt -\lbb i \int_0^1 h(X_{u,r,\ve}, \wt{y}_\ve^u)dr dt
        + f(u_\ve)  \wt{y}^u_\ve dt
        + e^{-W} R (\ve, \vf) dt,
\end{align}
where   $h(X_{u,r,\ve}, \wt{y}^u_\ve)$ and $f(u^\ve)$ are similar to those arising in \eqref{equa-z},
with  $X^u$, $z$, $u$ replaced by $X_{u,r,\ve}$, $\wt{y}^u_\ve$ and $u^\ve$ respectively.

Choose the Strichartz pair $(p,q)=(\a+1, \frac{4(\a+1)}{d(\a-1)})$. Then, by
Strichartz's estimates, H\"older's inequality and Minkowski's inequality, we have
\begin{align*}
   &\|\wt{y}^u_\ve\|_{C([0,t];L^2)} + \|\wt{y}^u_\ve\|_{L^q(0,t; L^p)}   \\
  \leq& C   \int_0^1 \|h(X_{u,r,\ve}, \wt{y}^u_\ve)\|_{L^{q'}(0,t,L^{p'})} dr
   + C t  \|f(u_\ve) \|_{L^\9(0,T;L^\9)} \|\wt{y}^u_\ve\|_{C([0,t];L^2)} \\
   &+ C  \|R (\ve, \vf)\|_{L^1(0,t; L^2) + L^{q'}(0,t; L^{p'})}.
\end{align*}
Note that,
\begin{align} \label{esti-hi-lpq}
        &\|h(X_{u,r,\ve}, \wt{y}^u_\ve)\|_{L^{q'}(0,t,L^{p'})}  \nonumber  \\
       \leq&  \a t^\theta \|\wt{y}^u_\ve\|_{L^q(0,t; L^p)} \|X_{u,r,\ve}\|^{\a-1}_{ L^{q}(0,t; L^{p})}  \nonumber  \\
       \leq&   \a (1\vee 2^{\a-1})  t^\theta \|\wt{y}^u_\ve\|_{L^q(0,t; L^p)} (\|X^{u_\ve} \|^{\a-1}_{L^q(0,T; L^p)} + \|X^u \|^{\a-1}_{L^q(0,T; L^p)} ),
\end{align}
where $\theta = 1 - d(\a-1)/4 >0$. Then,
\begin{align*}
       \|\wt{y}^u_\ve\|_{C([0,t];L^2)} + \|\wt{y}^u_\ve\|_{L^q(0,t; L^p)}
     \leq&  CD_3(T) (t^\theta + t) (\|\wt{y}^u_\ve\|_{C([0,t]; L^2)} + \|\wt{y}^u_\ve\|_{L^q(0,t; L^p)})  \\
         &  +C  \|R (\ve, \vf)\|_{L^1(0,T; L^2) + L^{q'}(0,T; L^{p'})}.
\end{align*}
where $D_3(T) = \a 2^{\a+1} \sup_{\ve \in[0,1]} (\|X^{u_\ve}\|^{\a-1}_{L^\9(\Omega; L^q(0,T; L^p))}
+ \|f(u_\ve)\|_{L^\9(\Omega; L^\9(0,T; L^\9))})$.
Using similar iterating arguments as in the proof of \eqref{bdd-forward-y}, we obtain
\begin{align*}
    \sup\limits_{u\in \calu_{ad}}(\|\wt{y}^u_\ve\|_{C([0,T];L^2)} + \|\wt{y}^u_\ve\|_{L^q(0,T; L^p)})
    \leq C(T)\|  R (\ve, \vf)\|_{L^1(0,T; L^2) + L^{q'}(0,T; L^{p'})}
\end{align*}
with $C(T) \in L^\9(\Omega)$.

Thus, by H\"older's inequality,
\begin{align*}
     & \bbe \|\wt{y}^u_\ve\|^2_{C([0,T];L^2)} + \bbe \|\wt{y}^u_\ve\|^2_{L^q(0,T; L^p)}   \\
     \leq& C(T)  \bbe  \|R (\ve, \vf)\|^2_{L^1(0,T; L^2) + L^{q'}(0,T; L^{p'})} \\
     \leq& C(T) \big( \ve^2 D_U^2 T^2 \|V\|^2_{L^\9(0,T;\bbr^m)}  \bbe \|\vf\|^2_{C([0,T]; L^2)}
                   +\sum\limits_{j=1}^2 \bbe \|R_j(\ve)\vf\|^2_{L^{q'}(0,T; L^{p'})}  \big).
\end{align*}

Therefore, in order to prove \eqref{asy-X-eta}, we only need to show that
\begin{align} \label{conv-R1}
   \bbe  \|R_j(\ve) \vf\|^2_{L^{q'}(0,T;L^{p'})} \to 0,\ \ as\ \ve\to 0,\ j=1,2.
\end{align}

Below, we prove \eqref{conv-R1} only for $R_1(\ve)$, but the argument applies as well to $R_2(\ve)$.
As in the proof of \eqref{conv-y-Xn*} we get
\begin{align} \label{conv-Xue-Xu}
     \|X^{u_\ve} - X^u\|_{C([0,T]; L^2)} + \|X^{u_\ve} - X^u\|_{L^q(0,T; L^{p})}
    \to 0,\ \ as \ \ve\to 0,\ \ \bbp-a.s.
\end{align}

Note that
\begin{align} \label{diff-h1}
  h_1(X_{u,r,\ve}) - h_1(X^u)
  =& \int_0^1 \partial_z h_1(X^u + r' (X_{u,r,\ve}-X^u))dr' (X_{u,r,\ve} - X^u) \nonumber  \\
   & + \int_0^1 \partial_{\ol{z}} h_1(X^u + r' (X_{u,r,\ve}-X^u))dr' (\ol{X_{u,r,\ve}} - \ol{X^u}).
\end{align}
Since $|\partial_z h_1(z)| + |\partial_{\ol{z}} h_1(z)| \leq C |z|^{\a-2}$ for $z\in \mathbb{C}$, using
the Minkowski inequality and the H\"older inequality we get that $\bbp$-a.s. for each $r\in [0,1]$,
\begin{align} \label{esti-Xue-Xu}
   &\| h_1(X_{u,r,\ve}) - h_1(X^u)\|_{L^{\frac{q}{\a-1}}(0,T; L^{\frac{p}{\a-1}})} \nonumber \\
   \leq& C \int_0^1  \|X^u + r' (X_{u,r,\ve}-X^u)\|^{\a-2}_{L^q(0,T; L^p)} dr' \|X_{u,r,\ve} - X^u\|_{L^q(0,T; L^p)} \nonumber \\
   \leq& C \sup\limits_{\ve\in[0,1]} \|X^{u_\ve}\|^{\a-2}_{L^q(0,T; L^p)} \|X^{u_\ve} - X^u\| _{L^q(0,T; L^p)}
   \to 0,\ \ as\ \ve \to 0,
\end{align}
where we also used $\a\geq 2$ and the last step is due to \eqref{conv-Xue-Xu}.

Thus, using the H\"older inequality combined with the Minkowski inequality and the bounded dominated convergence theorem we obtain that
\begin{align}\label{conv-h1-lpq}
    \|R_1(\ve)\vf\|_{L^{q'}(0,T; L^{p'})}
    \leq& T^\theta \|\vf\|_{L^q(0,T; L^p)} \int_0^1 \| h_1(X_{u,r,\ve}) - h_1(X^u)\|_{L^{\frac{q}{\a-1}}(0,T; L^{\frac{p}{\a-1}})}  dr \nonumber  \\
    \to& 0,\ \ as\ \ve \to 0, \ \bbp-a.s.
\end{align}

Moreover, taking into account \eqref{esti-Xue-Xu}, \eqref{Cons-bdd-Xn-Lpq} and  \eqref{bdd-forward-y}, we have
\begin{align*}
      \|R_1(\ve) \vf \|_{L^{q'}(0,T;L^{p'})}
     \leq   C T^{ \theta}\sup\limits_{\ve\in[0,1]}\|X^{u_\ve}\|^{ \a-1 }_{L^q(0,T; L^p)}   \|\vf\| _{L^q(0,T; L^{p})} \in L^\9(\Omega),
\end{align*}
which along with  \eqref{conv-h1-lpq} and the bounded dominated convergence theorem
yields \eqref{conv-R1} for $j=1$.  Therefore, the proof is complete.  \hfill $\square $\\

We shall prove now the following result.

\begin{proposition} \label{WP-Backward}
$(i)$ Assume Hypothesis $(H2)$ and that
$\bbx_T \in L^{2+\nu}(\Omega; L^2)$,  $\bbx_1\in L^{2+\nu}(\Omega; L^2(0,T; L^2))$ for
some small $\nu\in(0,1)$.

Then,
there exists a unique $(\mathcal{F}_t)$-adapted solution $(Y^u,Z^u)$
to \eqref{back-equa} corresponding to $u\in \calu_{ad}$,
satisfying for  any Stirchartz pair $(p,q)$,
\begin{align} \label{bdd-Yn}
     \sup\limits_{u\in \calu_{ad}}  (\|Y^u\|_{L^{2+\nu}(\Omega; C([0,T]; L^2))} + \|Y^u\|_{L^{2+\nu}(\Omega; L^q(0,T; L^p))}) < \9,
\end{align}
and
\begin{align} \label{bdd-Zn}
     \sup\limits_{u\in \calu_{ad}}   \|Z_k^u\|_{L^{2+\nu}(\Omega; L^2(0,T; L^2))} < \9,\ \ 1\leq k\leq N.
\end{align}

$(ii)$ Assume in addition  that $\bbx_T \in L^{2+\nu}(\Omega; H^1)$, $\bbx_1\in L^{2+\nu}(\Omega; L^2(0,T; H^1))$
for some small $\nu \in (0,1)$.
Then,  for  any $\rho\in [2,2+\nu)$ and any Strichartz
pair $(p,q)$, we have
\begin{align} \label{bdd-Yn-Wpq}
     \sup\limits_{u\in \calu_{ad}}  (\|Y^u\|_{L^{\rho}(\Omega; C([0,T]; H^1))} + \|Y^u\|_{L^{\rho}(\Omega; L^q(0,T; W^{1,p}))}) < \9,
\end{align}
and
\begin{align} \label{bdd-Zn-Wpq}
     \sup\limits_{u\in \calu_{ad}}   \|Z_k^u\|_{L^{\rho}(\Omega; L^2(0,T; H^1))} < \9,\ \ 1\leq k\leq N.
\end{align}

\end{proposition}

As mentioned in Section \ref{PROBLEM-RESULT}, the main difficulty in the analysis of backward stochastic equation
mainly comes from the singular term $\lbb i h_2(X^u)\ol{Y}$. The proof of Proposition \ref{WP-Backward}
follows the following steps.

First, we will consider  the truncated approximating equation \eqref{appro-back-equa} and
introduce the dual equation \eqref{forward-equa-Psi} below, which are related together by the formula \eqref{dual-Psi-Y}.
Then, the uniform estimate \eqref{psi-Psi} of dual solutions imply, via duality arguments, those of the approximating solutions $\{Y_n\}$
(see \eqref{esti-Yn-Lp} below), which in turn imply the uniform estimates of $\{Z_n\}$ (see \eqref{esti-Zn-L2} below).
Consequently, in view of the linear structure of \eqref{appro-back-equa}, one can pass to the limit and obtain the well-posedness of problem \eqref{back-equa}
as well as the estimates \eqref{bdd-Yn} and \eqref{bdd-Zn}.

Analogous arguments are applicable to prove estimates \eqref{bdd-Yn-Wpq} and \eqref{bdd-Zn-Wpq} in Sobolev spaces,
which requires the condition $\a\geq 2$ and the  integrability conditions on $\bbx_T$ and $\bbx_1$ in Sobolev spaces. \\

{\it \bf Proof. } $(i)$. Let $g$ be a radial smooth cut-off function  such that $g=1$ on $B_1(\bbr)$, and $g=0$ on $B^c_2(\bbr)$. For $j=1,2$,
set $h_{j,n}(X^u):= g(\frac{|X^u|}{n}) h_j(X^u)$.
Note that $ |h_{1,n}(X^u)| +  |h_{2,n}(X^u)| \leq \a 2^{\a-1} |g|_{L^\9} n^{\a-1}$.

Consider the approximating backward stochastic equation
\begin{align} \label{appro-back-equa}
 &d Y_n= -i\Delta Y_n\,dt - \lbb i h_{1,n}(X^u)Y_n dt +\lbb i  h_{2,n}(X^u)  \ol{Y_n} dt + \mu Y_n dt - iV_0Y_n dt \nonumber  \\
 &\qquad \quad        - i u\cdot V Y_n dt + \g_1 (X^u - \bbx_1) dt  - \sum\limits_{k=1}^N \ol{\mu_k} e_k Z_{k,n} dt +  \sum\limits_{k=1}^N Z_{k,n} d\beta_k(t),\nonumber    \\
 & Y_n(T) = -(X^u(T)-\bbx_T).
\end{align}
By standard theory for stochastic backward infinite dimensional equations (see, e.g., \cite{9a}, \cite{HP91}), it follows that
there exists a
unique $(\mathcal{F}_t)$-adapted solution $(Y_n,Z_n) \in  L^2(\Omega; C([0,T]; L^2))  \times (L^2_{ad}(0,T; L^2(\Omega; L^2)))^N$ to \eqref{appro-back-equa}.

In order to pass to the limit $n\to \9$,   we are going to obtain uniform estimates of $Y_n$ in the space
$\caly_{\rho_\nu}:= L^{\rho_\nu}(\Omega; L^\9(0,T; L^2)) \cap L^{\rho_\nu} (\Omega; L^q(0,T; L^p))$ with $\rho_\nu := 2+\nu$ and
$(p,q)=(\a+1, \frac{4(\a+1)}{d(\a-1)})$.

To this purpose, for each $n\geq 1$, define the functional $\Lambda_n$ on the space
$L^\9(\Omega\times (0,T) \times \bbr^d)$,
\begin{align} \label{Oper-Phi}
   \Lambda_n(\Psi) := \bbe& Re \<X^u(T)- \bbx_T, \psi_n(T)\>_2
                      +\g_1 \bbe \int_0^{T} Re\<X^u(t)-\bbx_1(t), \psi_n(t)\>_2dt,
\end{align}
where $\Psi \in L^\9(\Omega\times (0,T) \times \bbr^d)$, $X^u$ is the solution to \eqref{equa-x}, and $\psi_n$  satisfies
\begin{align}  \label{forward-equa-Psi}
  &d \psi_n = -i \Delta \psi_n dt - \lbb i h_{1,n}(X^u) \psi_n dt - \lbb i h_{2,n}(X^u) \ol{\psi_n} dt -\mu \psi_n dt  \nonumber \\
  &\qquad \quad - i V_0 \psi_n dt  -i u\cdot V \psi_n dt  - \Psi dt + \psi_n d W(t), \nonumber \\
  &\psi_n(0)= 0.
\end{align}
(Note that,  \eqref{forward-equa-Psi} is similar  to \eqref{forward-equa}
but with $\Psi$ in place of $i\wt{u}\cdot V X^u$.)

Since $|h_{j,n}(X^u)|\leq \a |g|_{L^\9} |X^u|^{\a-1}$, $j=1,2$, arguing as in the proof of Lemma \ref{WP-Forward}
we infer that there exists a unique strong $L^2$-solution $\psi_n$ to \eqref{forward-equa-Psi} on $[0,T]$, satisfying
\begin{align} \label{psi-Psi*}
     \sup\limits_{n} (\|\psi_n\|_{  C([0,T]; L^2)} + \|\psi_n\|_{ L^q(0,T; L^p )})
   \leq  C(T)  \|\Psi\|_{L^1(0,T; L^2) + L^{q'}(0,T;L^{p'})},
\end{align}
where $C(T) \in L^\9(\Omega)$ is independent of $n$ and $\Psi$.
It follows also that for $\rho >1$,
\begin{align} \label{psi-Psi}
   \sup\limits_{n\geq 1} \|\psi_n\|_{\caly_{\rho'}}
   \leq C(\rho,T)  \|\Psi\|_{\mathcal{Y}'_{\rho'}},
\end{align}
where $C(\rho, T)$ is independent of $n$ and $\Psi$ and $\caly'_{\rho'}:= L^{\rho'}(\Omega; L^1(0,T; L^2))  + L^{\rho'}(\Omega; L^{q'}(0,T;L^{p'}))$.

Moreover, by It\^o's formula, for every $n\geq 1$ and $\Psi\in L^\9(\Omega\times (0,T) \times \bbr^d)$, we have
\begin{align} \label{dual-Psi-Y}
    \Lambda_n(\Psi)= \bbe \int_0^{T} Re \<\Psi, Y_n\>_2 dt.
\end{align}

Thus, by the conservation of $|X^u(t)|_2$ and by estimates \eqref{Oper-Phi}, \eqref{psi-Psi}, we have
\begin{align} \label{dual-esti}
   |\Lambda_n(\Psi)|
   \leq&  \g_1 \|X^u-\bbx_1\|_{L^{\rho_\nu}(\Omega; L^2(0,T; L^2))} \|\psi_n\|_{L^{\rho_\nu'}(\Omega; L^2(0,T; L^2))} \nonumber  \\
       & + \|X^u(T) - \bbx_T\|_{L^{\rho_\nu}(\Omega; L^2)}   \|\psi_n(T)\|_{L^{\rho_\nu'}(\Omega; L^2)}  \nonumber  \\
   \leq&  C \|\Psi\|_{\mathcal{Y}'_{\rho_\nu'}},
\end{align}
where $C$ is independent of $n$. Since $L^\9( \Omega\times (0,T) \times \bbr^d )$ is dense in $\caly'_{\rho_\nu'}$
and $\mathcal{Y}_{\rho_{\nu}}$ is the dual space of $\caly'_{\rho_\nu'}$,
it follows by \eqref{dual-Psi-Y}, \eqref{dual-esti} that
\begin{align}  \label{esti-Yn-Lp}
       \sup\limits_{n\geq 1}\|  Y_n\|_{\mathcal{Y}_{\rho_\nu}} <\9.
\end{align}

Hence,  there exists  $\wt{Y} \in \mathcal{Y}_{\rho_\nu}$, such that along a subsequence of
$\{n\}\to \9$ (still denoted by $\{n\}$),
\begin{align} \label{conv-Yn-Lp}
    Y_n  \to \wt{Y},\ \ weak-star\ in\    \mathcal{Y}_{\rho_\nu}.
\end{align}

Note that, for each $j=1,2$, $h_{j,n}(X^u) \to h_j(X^u)$,  $d\bbp \otimes dt \otimes d\xi-a.e.$,
and $\sup_{n\geq 1}|h_{j,n}(X^u)| \leq \a |g|_{L^\9} |X^u|^{\a-1} \in L^{2\rho_\nu'}(\Omega; L^{\frac{q}{\a-1}}(0,T;L^{\frac{p}{\a-1}}))$,
by the dominated convergence theorem,
\begin{align} \label{conv-hin}
    h_{j,n}(X^u) \to h_j(X^u), \ in\ L^{2\rho_\nu'}(\Omega; L^{\frac{q}{\a-1}}(0,T;L^{\frac{p}{\a-1}})), \ \ as\ n\to \9.
\end{align}
Since $\frac{1}{(2\rho'_\nu)'} = \frac{1}{\rho_\nu} + \frac{1}{2\rho_\nu}$, from H\"older's inequality, \eqref{conv-Yn-Lp} and \eqref{conv-hin}  it follows  that
\begin{align} \label{conv-hin-Yn}
    h_{1,n} (X^u) Y_n \to  h_1(X^u) \wt{Y},\  h_{2,n} (X^u) \ol{Y_n} \to  h_2(X^u) \ol{\wt{Y}},
\end{align}
weakly in $L^{(2\rho_\nu')'}(\Omega; L^{\frac{q}{\a}}(0,T;L^{p'}))$.

Moreover, we claim that for $1\leq k\leq N$,
\begin{align} \label{esti-Zn-L2}
   \sup\limits_{n\geq 1} \|Z_{k,n}\|_{L^{\rho_\nu}(\Omega; L^2(0,T; L^2))} \leq C <\9.
\end{align}
In particular, for each $1\leq k\leq N$,  there exists $Z^u_k\in L^{\rho_\nu}(\Omega; L^2(0,T; L^2))$  such that
(selecting a further subsequence if necessary)
\begin{align} \label{conv-Zn-L2}
    Z_{k,n} \to Z^u_k,\ \ weakly\ in\ L^2(\Omega; L^2(0,T; L^2)).
\end{align}
Since $v \mapsto \int_\cdot^T  v d\beta_k(s)$ is a bounded linear operator in  $L^2(\Omega; L^2(0,T; L^2))$, it follows that
for each $1\leq k\leq N$,
\begin{align} \label{conv-Zn-L2*}
   \int_\cdot^T Z_{k,n} d\beta_k(s) \to   \int_\cdot^T Z^u_{k} d\beta_k(s),\ \ weakly\ in\ L^2(\Omega; L^2(0,T; L^2)).
\end{align}

To prove \eqref{esti-Zn-L2}, we apply It\^o's formula to \eqref{appro-back-equa} to get that for $\eta>0$,
\begin{align} \label{ito-pn}
      &e^{\eta t}|Y_n(t)|_2^2  \nonumber \\
     =& e^{\eta T} |X^{u}(T) - \bbx_T|_2^2
        - \eta \int_t^T e^{\eta s} |Y_n|_2^2 ds   \nonumber \\
      & + 2\lbb \int_{t}^{T}  Im \int e^{\eta s} h_{2,n}(X^u) \ol{Y_n^2} d\xi ds
        -2 \int_{t }^{T}  \int e^{\eta s} \mu |Y_n|^2 d\xi ds \nonumber \\
       & -2 \g_1 \int_{t }^{T} Re  \int e^{\eta s} Y_n (\ol{X^{u}} - \ol{\bbx_1}) d\xi  ds
       +2  \sum\limits_{k=1}^N \int_{t }^{T}   Re  \int  e^{\eta s}  \mu_k e_k Y_n\ol{Z_{k,n}} d\xi ds \nonumber  \\
      &  - \sum\limits_{k=1}^N \int_{t }^{T} e^{\eta s} |Z_{k,n}|_2^2 ds
         -2  \sum\limits_{k=1}^N\int_{t }^{T} Re \int e^{\eta s} Y_n \ol{Z_{k,n}} d\xi  d\beta_k(s)     \nonumber \\
      =:& e^{\eta T} |X^{u}(T) -  \bbx_T|_2^2 - \eta \int_t^T e^{\eta s} |Y_n|_2^2 ds  + \sum\limits_{j=1}^6 K_j(t ).
\end{align}

Note that, since $|h_{2,n}(X^u)| \leq \a |g|_{L^\9} |X^u|^{\a-1}$, by H\"older's inequality,
\begin{align} \label{esti-K1}
       |K_1(t)|
       \leq  \a |g|_{L^\9} e^{\eta T}  T^{\theta} \|X^{u}\|^{\a-1}_{L^q(0,T; L^{p})} \|Y_n\|^2_{L^q(0,T; L^p)},
\end{align}
where $\theta  = 1 - d(\a-1) / 4 \in (0,1)$.

Moreover, using $ab\leq ca^2 + c^{-1}b^2$, $c>0$, we get
\begin{align}  \label{esti-K24}
     \sum\limits_{j=2}^4 |K_j(t)|
     \leq&  (2|\mu|_\9 + 2 \g_1 + 8\sum\limits_{k=1}^N|\mu_ke_k|^2_{L^\9}) \int_{t }^{T} e^{\eta s} |Y_n|_2^2 ds  \nonumber  \\
         & +   2 \g_1 \int_{t }^{T} e^{\eta s} |X^{u} - \bbx_1|_2^2 ds
           +  \frac 12 \int_{t}^{T} e^{\eta s} |Z_{k,n}|_2^2 ds.
\end{align}

Thus, choosing  $\eta > 2|\mu|_\9 + 2 \g_1 + 8\sum_{k=1}^N|\mu_ke_k|^2_{L^\9}$, it follows  that for any $t\in[0,T]$,
\begin{align} \label{esti-YZ}
    e^{\eta t}|Y_n(t)|_2^2 + \frac 12 \sum\limits_{k=1}^N\int_{t}^{T} e^{\eta s} |Z_{k,n}|_2^2 ds
   \leq  V_{T,n} -2  \sum\limits_{k=1}^N\int_{t }^{T} Re \int e^{\eta s} Y_n \ol{Z_{k,n}} d\xi  d\beta_k(s),
\end{align}
where
\begin{align} \label{VT}
   V_{T,n} : =& e^{\eta T} |X^{u}(T) - \bbx_T|_2^2
         + 2\g_1  \int_{0}^{T} e^{\eta s} |X^{u} - \bbx_1|_2^2 ds \nonumber \\
       & + \a |g|_{L^\9} e^{\eta T}  T^{\theta} \|X^{u}\|^{\a-1}_{L^q(0,T; L^{p})} \|Y_n\|^2_{L^q(0,T; L^p)}.
\end{align}
This yields
\begin{align} \label{esti-Z}
       & \sum\limits_{k=1}^N \bbe \(\int_{0}^{T} e^{\eta s} |Z_{k,n}|_2^2 ds\)^{\frac{\rho_\nu}{2}} \nonumber \\
   \leq& C (\rho_\nu) \( \bbe V^{\frac{\rho_\nu}{2}}_{T,n} + \bbe \sup\limits_{t\in[0,T]}
         \bigg|\sum\limits_{k=1}^N\int_{0 }^{t} Re \int e^{\eta s} Y_n \ol{Z_{k,n}} d\xi  d\beta_k(s)\bigg|^{\frac{\rho_\nu}{2}} \).
\end{align}

Note that, $\bbe V^{\rho_\nu/2}_{T,n} <\9$, due to \eqref{Cons-bdd-Xn-Lpq}, \eqref{esti-Yn-Lp} and to the integrability conditions on $\bbx_T$ and $\bbx_1$. Moreover,
by the Burkholder-Davis-Gundy inequality, the
second term in the right hand side of \eqref{esti-Z} is bounded by
\begin{align*}
     & C (\rho_\nu) \bbe \(\int_{0 }^{T}  \sum\limits_{k=1}^N \bigg|\int e^{\eta s} Y_n \ol{Z_{k,n}} d\xi \bigg|^2 ds\)^{\frac{\rho_\nu}{4}} \\
   \leq&  C(\rho_\nu) \sum\limits_{k=1}^N \bbe  \sup\limits_{t\in[0,T]}  e^{\frac 14 \rho_\nu \eta t} |Y_n(t)|_2^{\frac{\rho_\nu}{2}}
         \( \int_0^T e^{\eta s}|Z_{k,n}|_2^2 ds \)^{\frac {\rho_\nu}{4}} \\
   \leq&  C(\rho_\nu, N) e^{\frac 12 \rho_\nu\eta T}  \bbe |Y_n|_{C([0,T]; L^2)}^{ \rho_\nu}
         + \frac 12 \sum\limits_{k=1}^N \bbe \(\int_{0}^{T} e^{\eta s} |Z_{k,n}|_2^2 ds\)^{\frac{\rho_\nu}{2}}.
\end{align*}
Plugging this into \eqref{esti-Z} we get
\begin{align} \label{esti-Z-Y}
   & \sup\limits_{n\geq 1}\frac 12 \sum\limits_{k=1}^N \bbe \(\int_{0}^{T} e^{\eta s} |Z_{k,n}|_2^2 ds\)^{\frac{\rho_\nu}{2}}  \nonumber \\
    \leq&  C(\rho_\nu) \sup\limits_{n\geq 1}\bbe V^{\frac{\rho_\nu}{2}}_{T,n}
         +  C(\rho_\nu, N) e^{\frac 12 \rho_\nu \eta T} \sup\limits_{n\geq 1}\bbe|Y_n|_{C([0,T]; L^2)}^{\rho_\nu},
\end{align}
which by \eqref{esti-Yn-Lp} implies \eqref{esti-Zn-L2}, as claimed.

Now, set
\begin{align*}
   Y^u(\cdot):=&   -(X^u(T) - \bbx_T) - \int_{\cdot }^{T}  \bigg(-i \Delta \wt{Y} - \lbb i h_1(X^u) \wt{Y}
            +\lbb i  h_2(X^u)  \ol{\wt{Y}} \\
           &\qquad + \mu \wt{Y}
         - iV_0 \wt{Y}  - i u\cdot V \wt{Y}   + \g_1 (X^u - \bbx_1)   - \sum\limits_{k=1}^N \ol{\mu_k} e_k Z^u_{k} \bigg) ds \\
           &-  \sum\limits_{k=1}^N \int_{\cdot }^{T} Z^u_{k} d\beta_k(s) .
\end{align*}
By virtue of \eqref{conv-Yn-Lp}, \eqref{conv-hin-Yn}, \eqref{conv-Zn-L2} and \eqref{conv-Zn-L2*},
we may pass to the limit in \eqref{appro-back-equa} and obtain that for any $v\in H^2$, $f \in L^\9(\Omega\times(0,T))$,
\begin{align}\label{back-equa-n.1}
    \bbe \int_0^T {}_{H^{-2}}\<\wt{Y}(t), f(t)v\>_{H^2} dt
    = \bbe \int_0^T {}_{H^{-2}}\<Y^u(t), f(t)v\>_{H^2} dt,
\end{align}
which implies that $Y^u=\wt{Y}$, in $H^{-2}$,
$d\bbp \otimes dt$-a.e. Since $Y^u$ is continuous in $H^{-2}$, $\bbp$-a.s., we can find a  null set $N'$,
such that for any $\omega\not \in N'$, $(Y^u(\omega), Z^u(\omega))$ solves \eqref{back-equa} in $H^{-2}$
for all $t\in [0,T]$, which proves the existence of solution to \eqref{back-equa}.
Estimates \eqref{bdd-Yn} and \eqref{bdd-Zn}  follow immediately by
\eqref{esti-Yn-Lp} and \eqref{esti-Zn-L2} respectively. Moreover, as in the proof of \cite[Lemma 4.3]{1},
we have for $|Y^u(t)-Y^u(s)|_2^2$ an It\^o formula
similar to \eqref{ito-pn}, which implies $t\to Y^u(t)$ is $L^2$ continuous.

The uniqueness can also be proved by the duality arguments. Indeed, let $(Y^u_j, Z_j^u)$, $j=1,2$, be any two solutions to \eqref{back-equa}.
For any $\Psi\in L^\9(\Omega \times (0,T) \times \bbr^d)$, let
$\psi$ be the unique solution to \eqref{forward-equa-Psi} but with $h_{j,n}(X^u)$ replaced by $h_j(X^u)$, $j=1,2$. Define $\Lambda(\Psi)$
similarly as in \eqref{Oper-Phi} with $\psi_n$ replaced by $\psi$. Then, similarly to \eqref{dual-Psi-Y},
$\Lambda(\Psi) = \bbe\int_0^T Re \<\Psi, Y^u_j\>_2dt$, $j=1,2$. It follows that $Y^u_1 = Y^u_2$ by
the arbitrariness of $\Psi $, and so $Z^u_1 = Z^u_2$ by the
estimate similar to \eqref{esti-Z-Y}. \\

$(ii)$. Fix $1\leq j \leq N$.  Consider the approximating equation of $(\partial_j Y^u, \partial_j Z^u)$ below ($\partial_j := \frac{\partial}{\partial \xi_j}$),
\begin{align} \label{appro-back-equa-h1}
& d Y'_n=  -i\Delta  Y'_n dt + G_n (Y'_n) dt
                     - \sum\limits_{k=1}^N \ol{\mu}_k e_k Z'_{k,n} dt + \g_1 \partial_j (X-\bbx_1) dt  \nonumber  \\
& \qquad \qquad \  + F_n(X^u, Y^u, Z^u) dt +  \sum\limits_{k=1}^N  Z'_{k,n} d\beta_k(t),  \nonumber  \\
& Y'_n=  - (\partial_j X^u(T) - \partial_j \bbx_T),
\end{align}
where $X^u$ and $(Y^u, Z^u)$ are the solutions to \eqref{equa-x} and \eqref{back-equa} respectively,
\begin{align*}
    G_n (Y'_n )
    := - \lbb i h_{1,n}(X^u)Y'_n  + \lbb i h_{2,n}(X^u) \ol{Y'_n }
       +( \mu - i V_0 - iu\cdot V) Y'_n ,
\end{align*}
\begin{align*}
   F_n(X^u, Y^u, Z^u) :=& -\lbb i h'_{1,n}(X^u)Y^u + \lbb i h'_{2,n}(X^u) \ol{Y^u}
                   +  \partial_j \mu Y^u   - i \partial_j V_0Y^u \\
                   & - i u\cdot \partial_j V Y^u
                     - \sum\limits_{k=1}^N \partial_j (\ol{\mu_k} e_k) Z^u_{k},
\end{align*}
$h_{k,n}(X^u)$ is as in \eqref{appro-back-equa} and $h'_{k,n}(X^u) = g(\frac{|X^u|+|\na X^u|}{n}) \partial_j(h_{k}(X^u))$, $k=1,2$. By truncation,
$|h_{k,n}(X^u)| +| h'_{k,n}(X^u)|  \leq C n^{\a-1}$, $k=1,2$, it follows
that there exists a unique $(\mathcal{F}_t)$-adapted solution $(Y'_n, Z'_n) \in L^2(\Omega; C([0,T]; L^2)) \times (L^2_{ad}(0,T; L^2(\Omega; L^2)))^N$
to \eqref{appro-back-equa-h1}.

For each $\Psi\in L^\9(\Omega \times (0,T) \times \bbr^d)$, let $\psi_n$ be the solution to \eqref{forward-equa-Psi}.
Similarly to \eqref{Oper-Phi},  set
\begin{align*}
   \wt{\Lambda}_{j,n}(\Psi) :=& \bbe Re \<\partial_j X^u(T)- \partial_j \bbx_T, \psi_n(T)\>_2 \\
                     & + \g_1 \bbe \int_0^{T} Re\< \partial_j X^u(t)- \partial_j \bbx_1(t), \psi_n(t))\>_2dt .
\end{align*}
By It\^o's formula, we have
\begin{align} \label{def-Lamj}
\wt{\Lambda}_{j,n}(\Psi) = \bbe \int_0^T Re\<\Psi, Y'_n\>_2 dt - \bbe \int_0^T Re \int F_n(X^u,Y^u, Z^u)\ol{\psi_n} d\xi dt.
\end{align}

Note that, since $|h'_{1,n}(X^u)| \leq C |X^u|^{\a-2} |\partial_j X^u|$ and $2\leq \a < 1+ 4/d$,
by H\"older's inequality, we have for $(p,q)=(\a+1, \frac{4(\a+1)}{d(\a-1)})$,
\begin{align} \label{restrict}
    & \bigg|  \int_0^T Re \int (-\lbb i)  h'_{1,n}(X^u) Y^u \ol{\psi_n} d\xi dt \bigg| \nonumber \\
    \leq&   \|h'_{1,n}(X^u)Y^u \|_{L^{q'}(0,T; L^{p'})} \|\psi_n\|_{L^q(0,T;L^p)}  \nonumber \\
    \leq& C_\a   T^\theta \|X^u\|^{\a-2}_{L^q(0,T;L^p)} \|\partial_j X^u\|_{L^q(0,T;L^p)}\|Y^u \|_{L^q(0,T;L^p)} \|\psi_n\|_{L^q(0,T;L^p)} \nonumber  \\
    \leq& C_\a   T^\theta \|X^u\|^{\a-1}_{L^q(0,T;W^{1,p})} \|Y^u\|_{L^q(0,T;L^p)}  \|\psi_n\|_{L^q(0,T;L^p)},
\end{align}
where $\theta = 1- d(\a-1)/4\in (0,1)$. Hence, for any $\rho \in [2,\rho_\nu)$,
\begin{align*}
   &\bigg| \bbe \int_0^T Re \int (-\lbb i)  h'_{1,n}(X^u) Y^u \ol{\psi_n} d\xi dt \bigg|\nonumber \\
   \leq& C(T) \bbe (\|X^u\|^{\a-1}_{L^q(0,T;W^{1,p})} \|Y^u\|_{L^q(0,T;L^p)}  \|\psi_n\|_{L^q(0,T;L^p)}) \nonumber \\
   \leq& C(T) \|X^u\|^{\a-1}_{L^\eta(\Omega; L^q(0,T;W^{1,p}))}  \|Y^u\|_{L^{\rho_\nu}(\Omega; L^q(0,T;L^p))}    \|\psi_n\|_{L^{\rho'}(\Omega;L^q(0,T;L^p))}),
\end{align*}
where $\eta$ satisfies $\frac{1}{(\a-1)\eta} = \frac{1}{\rho'_\nu} - \frac{1}{\rho'}>0$.
Similar arguments apply to the term involving $\lbb i h'_{2,n}(X^u) \ol{Y^u}$.
Moreover, the other terms in the integration
$\bbe \int_0^T Re \int F_n(X^u,Y^u,Z^u )\ol{\psi_n} d\xi dt$ are bounded by
\begin{align*}
   C \(\|Y^u \|_{L^\rho(\Omega; L^\9(0,T;L^2))}  + \sum\limits_{k=1}^N\|Z^u_k\|_{L^\rho(\Omega; L^2(0,T;L^2))}\) \|\psi_n\|_{L^{\rho'}(\Omega; L^\9(0,T;L^2))}.
\end{align*}
Plugging the estimates above into \eqref{def-Lamj} and using \eqref{psi-Psi}, \eqref{Cons-bdd-Xn-Lpq}, \eqref{bdd-Yn} and \eqref{bdd-Zn}
we obtain for any  $\rho \in [2,\rho_\nu)$,
\begin{align*}
   \bigg|\bbe \int_0^T Re\<\Psi,Y'_n\>_2 dt\bigg|
   \leq&  |\wt{\Lambda}_{j,n}(\Psi)| +  \bigg|\bbe \int_0^T Re \int F_n(X^u,Y^u,Z^u)\ol{\psi_n} d\xi dt \bigg| \\
   \leq&   C \|\Psi\|_{\mathcal{Y}'_{\rho'}}
\end{align*}
with $C$ independent of $n$ and $\Psi$,
which implies that for any $\rho \in [2,\rho_\nu)$,
\begin{align} \label{esti-Yn'}
     \sup\limits_{n\geq 1} \| Y'_n \|_{\mathcal{Y}_{\rho}} \leq C.
\end{align}
Once we obtain \eqref{esti-Yn'}, using similar arguments as those below \eqref{esti-Yn-Lp},
we can prove the assertion $(ii)$. The details are omitted.   \hfill $\square$ \\

{\it \bf  Proof of Proposition \ref{l4.1}.} Using \eqref{asy-X-y} in Lemma \ref{WP-Forward} we have
\begin{align}\label{e4.9}
   & \lim_{\vp\to0}\frac1\vp\,(\Phi(u+\vp\wt u)-\Phi(u))   \nonumber \\
 =& 2\E \bigg(\!{\rm Re}\<X^u(T)-\bbx_T,\vf^{u,\wt{u}}(T)\>_2 \nonumber \\
  &\qquad  +\g_1\E\int^{T}_0 {\rm Re}\<X^u(t)-\bbx_1(t),\vf^{u,\wt{u}}(t)\>_2dt
     +\g_2\int^{T}_0 u\cdot{\wt u}\,dt \bigg).
\end{align}
Then, similarly to \eqref{dual-Psi-Y}, by \eqref{back-equa} and \eqref{forward-equa}, we obtain via It\^o's formula,
\begin{align*}
 \E\,{\rm Re}\<X^u(T)-\bbx_T,\vf^{u,\wt{u}}(T)\>_2
+&\g_1\E\dd\int^{T}_0{\rm Re}\<X^u(t)-\bbx_1(t),\vf^{u,\wt{u}}(t)\>_2  dt\nonumber  \\
&= -\E\,{\rm Im}\int^{T}_0\int_{\rr^d}\wt u\cdot V\,X^u  \ol{Y^u} d\xi dt.
\end{align*}
Combining these formulas we get  \eqref{e4.1}  as claimed.  \hfill $\square$

\section{Proof of Theorem \ref{t2.6}.} \label{PROOF-THM2}

As in the proof of Lemma \ref{Lem-conv}, we note that $\Phi$ is  continuous on the metric space $\calu_{ad}$ endowed with the distance
$d(u,v)= \|u-v\|=(\bbe \int_0^T |u(t)-v(t)|_m^2 dt)^{1/2}$.
Applying Ekeland's variational principle in $\calu_{ad}$ (see \cite[Theorem 1]{9}, or  \cite{8}),
for every $n\in\nn$  we get   $u_n\in\calu_{ad}$ such that
\begin{align}
 &\Phi(u_n)\le\Phi(u)+\dd\frac{1}{n}\,d(u_n,u),\ \ \forall  u\in\calu_{ad}.  \label{inf-phi-n}
\end{align}
In particular, it follows that
\begin{equation}\label{e4.10a}
u_n={\rm arg\,min}\left\{\Phi(u)+\frac{1}{n}\,\|u_n-u\|;\ u\in\calu_{ad}\right\}.
\end{equation}

We define the function $\wt{\Phi}: L^2_{ad}(0,T; \bbr^m) \to \ol{\bbr} = (-\9,+\9]$ by
\begin{align*}
    \wt{\Phi} = \Phi(u) + I_{\calu_{ad}}(u), \ \ \forall u\in L^2_{ad}(0,T; \bbr^m),
\end{align*}
where
\begin{align*}
    I_{\calu_{ad}} (u) = \left\{
                           \begin{array}{ll}
                             0, & \hbox{if $u\in \calu_{ad}$;} \\
                             +\9, & \hbox{otherwise.}
                           \end{array}
                         \right.
\end{align*}
The subdifferential $\partial \wt{\Phi} (u) \subset L^2_{ad}(0,T; \bbr^m)$ of $\wt{\Phi}$ at $u$ in the sense of
R.T. Rockafellar \cite{R79} is defined as the set of all $z\in L^2_{ad}(0,T; \bbr^m)$ such that the function
$v \to \wt{\Phi}(v)- \bbe \int_0^T v(t) z(t) dt$ has $u$ as a substationary point in the sense of \cite{R79}.

We have
\begin{align} \label{phi-eta}
    \partial \wt{\Phi} (u) \subset \eta(u) + \mathcal{N}_{\calu_{ad}} (u), \ \ \forall u \in \calu_{ad},
\end{align}
where $\eta (u)$ is defined by \eqref{e4.2}, $\mathcal{N}_{\calu_{ad}}(u)$ is the normal cone to
$\calu_{ad} $ at $u_n$, i.e.,
$$\mathcal{N}_{\calu_{ad}}(u_n)=\{v\in L^2_{ad}(0,T; \bbr^m) ;\  \<v,  u_n-{\wt v}\> \ge0,\ \ \ff\wt v\in \calu_{ad}\},$$
and $\<\ ,\ \>$ denotes the inner product of $ L^2_{ad}(0,T; \bbr^m)$.

To prove \eqref{phi-eta}, as mentioned in \cite[(2.4)]{R79}, for each $u\in \calu_{ad}$, one has
\begin{align*}
    \partial \wt{\Phi}(u)
    = \{ z\in  L^2_{ad}(0,T; \bbr^m): \wt{\Phi}^\uparrow (u, y) \geq \<y, z\>, \forall y \in  L^2_{ad}(0,T; \bbr^m) \},
\end{align*}
where $\wt{\Phi}^\uparrow (u,y)$ is the subderivative at $u$ with respect to $y$
\begin{align*}
    \wt{\Phi}^\uparrow (u, y)
    = \sup\limits_{V \subset \mathcal{N}(y)}
        \left[ \limsup\limits_{\substack{u'\to u,\a'\to \wt{\Phi}(u) \\ \a'\geq \wt{\Phi}(u'), t\to 0}}
              \inf\limits_{y'\in V}
              \left( \frac{\Phi(u'+ty')-\a'}{t}
        + \frac{I_{\calu_{ad}}(u'+ty')}{t} \right)  \right],
\end{align*}
and $\mathcal{N}(y)$ is the set of all neighborhoods of $y$. This yields
\begin{align*}
    \wt{\Phi}^\uparrow (u, y)
    = \lim\limits_{t\to 0} \frac{\Phi(u+ty)-\Phi(u)}{t} + I'_{\calu_{ad}}(u,y),
\end{align*}
where $I'_{\calu_{ad}}(u,y) =0$ if $y\in T_{\calu_{ad}}(u)$,
and $I'_{\calu_{ad}}(u,y) = \9$ if $y\not \in T_{\calu_{ad}}(u)$, $T_{\calu_{ad}}(u)$ is the (Clarke) tangent cone
to $\mathcal{U}_{ad}$ at $u$ defined in \cite{R79}.  Then, by Proposition \ref{l4.1},
for any $z\in \partial \wt{\Phi}(u)$, we have that $\<\eta(u), y\> \geq \<y,z\>$, $\forall y=v-u$, $v \in \calu_{ad}$.
Thus, $z\in \eta(u) + \mathcal{N}_{\calu_{ad}}(u)$. This implies \eqref{phi-eta}.

On the other hand, by Theorem $2$ in \cite{R79} we have
\begin{align} \label{subdiff-subset}
   \partial ( \wt{\Phi}(u) + \frac 1n \|u_n-u\| ) \subset \partial \wt{\Phi}(u) + \frac 1n \partial \|u_n-u\|.
\end{align}

Thus, by \eqref{e4.10a}-\eqref{subdiff-subset} we get
\begin{align*}
    0 \in& \partial (\wt{\Phi}(u)+ \frac 1n \|u_n-u\|) (u=u_n) \\
      \subset& \eta(u_n) + \frac 1n  (\partial \|u_n-u\|) (u=u_n) + \mathcal{N}_{\calu_{ad}}(u_n),
\end{align*}
which implies that there exist $\zeta_n \in \mathcal{N}_{\calu_{ad}}(u_n)$ and $\eta_n \in (\partial \|u_n-u\|) (u=u_n)$,
such that
\begin{align} \label{e4.11}
    \eta(u_n) + \zeta_n + \frac 1n \eta_n = 0.
\end{align}

We claim  that,
\begin{align} \label{equa-cone}
\mathcal{N}_{\calu_{ad}}(u_n)=\{v\in L^2_{ad}(0,T;\rr^m): v\in N_U(u_n),\ \ a.e.\ on\ (0,T) \times \Omega.\},
\end{align}
where $N_U(u_n)$ is the normal cone to $U\subset \rr^m$ at $u_n\in U,$ that is,
$N_U(u_n)=\{v\in\rr^m;\  v\cdot( u_n-{\wt v}) \ge0,\ \ \ff\wt v\in U\}.$

Indeed, for any $\eta \in \mathcal{N}_{\calu_{ad}}(u_n)$, we have
\begin{align}  \label{cone-eta}
     \bbe \int_0^T \eta \cdot (u_n - v) dt \geq 0, \ \ \forall\ v\in \calu_{ad}.
\end{align}
Since for each closed convex set $U$, $\forall \nu>0$, $(I+\nu N_U)^{-1}=P_U$,
where $P_U$ is the projection on $U$, there exists a unique $v \in \calu_{ad}$, such that
\begin{align} \label{equa-v}
    v + N_U(v) \ni u_n + \eta,\ \ a.e.\ on\ (0,T) \times \Omega,
\end{align}
i.e. $v = P_U (u_n + \eta)$,  a.e. on $(0,T) \times \Omega$.
Hence, there exists $\zeta_v \in N_U(v)$, such that $ v + \zeta_v = u_n + \eta$, $d \bbp \otimes dt $-a.e.
Then, by \eqref{cone-eta},
\begin{align*}
    0 \leq  \bbe \int_0^T (v-u_n + \zeta_v ) \cdot (u_n - v) dt
      = - \|u_n - v\|^2 + \bbe \int_0^T \zeta_v (u_n - v) dt.
\end{align*}
Since $d \bbp \otimes dt $-a.e., $\zeta_v \in N_U(v)$, $\zeta_v \cdot (v-u_n) \geq 0$, we get
\begin{align*}
      \|u_n - v\|^2 \leq \bbe \int_0^T \zeta_v\cdot (u_n -v) dt \leq 0.
\end{align*}
It follows that $u_n = v$, $d \bbp \otimes dt$-a.e., which yields by \eqref{equa-v} that $\eta \in N_U(u_n)$,  $d \bbp \otimes dt$-a.e.,
thereby implying
\begin{align*}
\mathcal{N}_{\calu_{ad}}(u_n) \subset \{v\in L^2_{ad}((0,T);\rr^m): v\in N_U(u_n),\ \ a.e.\ on\ (0,T) \times \Omega.\}.
\end{align*}
The inverse inclusion is obvious. Thus we obtain \eqref{equa-cone}, as claimed. \\

Now, by virtue of \eqref{e4.2}, we may rewrite \eqref{e4.11} as
\begin{equation}\label{e4.12}
\barr{r}
\dd u_n(t)+\frac1{2\gamma_2}\,\zeta^0_n(t)= \frac 1{\gamma_2}\,{\rm Im}\int_{\rr^d} V(\xi)  X_n(t,\xi) \ol{Y_n}(t,\xi)d\xi-\frac1{2\g_2 n}\,\eta_n(t)\vsp\mbox{ a.e. on }(0,T)\times\ooo,\earr\end{equation}
where $\zeta^0_n(t)\in N_U(u_n(t)),$  $X_n :=X^{u_n}$, and $(Y_n,Z_n)$ is the solution to
\eqref{back-equa} corresponding to $u_n$. We get by \eqref{e4.12} that
\begin{equation}\label{e4.13}
u_n(t)=P_U\(\dd\frac 1{\g_2}\,{\rm Im}\int V(\xi) X_n(t,\xi) \ol{Y_n}(t,\xi)d\xi-\frac1{2\g_2 n}\,\eta_n(t)\),
\end{equation}
with
\begin{equation}\label{e4.14}
\E\int^{T}_0|\eta_n(t)|^2_mdt=1.
\end{equation}

We claim that there exists a probability space $(\Omega^*, \mathcal{F}^*, \bbp^*)$,
$u^*_n, u^*\in \calu_{ad^*}$, $n\geq 1$, such that
the distributions of $u^*_n$ and $u_n$ coincide on $L^1(0,T; \bbr^m)$, and as $n\to \9$,
\begin{align} \label{un-u-L1}
    u^*_n \to u^*, \ \ in\ L^1(0,T; \bbr^m), \ \ \bbp^*-a.s.
\end{align}
Then, it follows from the boundedness of $\{u^*_n\}$ that
\begin{align*}
    u^*_n \to u^* , \ \ in\ L^2(0,T; \bbr^m), \ \ \bbp^*-a.s.
\end{align*}
Hence, similar arguments as in the proof of Theorem \ref{t2.5} imply that
\begin{align*}
    \Phi^*(u^*) = \lim\limits_{n\to \9} \Phi^*(u^*_n) =  \lim\limits_{n\to \9} \Phi(u_n) =  I,
\end{align*}
thereby yielding the sharp equality in \eqref{e2.10}.

It remains to prove \eqref{un-u-L1}. By virtue of Skorohod's representation theorem, we only need to show the tightness of
the distributions of $u_n$ in $L^1(0,T; \bbr^m)$, $n\geq 1$.
For this purpose, in view of Lemma \ref{Lem-tight-L1},
it suffices to prove that $\mu_n:=\bbp \circ u_n^{-1}$, $n\geq 1$,
satisfy \eqref{tight-space} and \eqref{tight-time}.

Indeed, \eqref{tight-space} follows immediately from the uniform boundedness of $\{u_n\}$. Regarding \eqref{tight-time},
by Markov's inequality, it suffices to show that there exists a positive exponent $b >0$ such that for any $\delta\in(0,1)$,
\begin{align} \label{L1-h-un}
   \limsup\limits_{n\to \9} \bbe \sup\limits_{0<h\leq \delta}  \int_0^{T-h}  | u_n(t+h) - u_n(t)|_m  dt  \leq C \delta^b.
\end{align}

To this end, since $P_U$ is Lipschitz, using \eqref{e4.13}, the Chauchy  inequality, \eqref{bdd-Xn-H1-Wpq} and \eqref{bdd-Yn-Wpq}  we get
\begin{align} \label{L1-h-un-esti}
  &\bbe \sup\limits_{0<h\leq \delta}  \int_0^{T-h} | u_n(t+h) - u_n(t) |_m dt \nonumber  \\
  \leq& \frac{1}{\g_2 n} T^{\frac 12} + \frac{1}{\g_2}  |V|_{L^{\9}} \bbe  \int_0^{T}  \sup\limits_{0<h\leq \delta}  \bigg( |X_n(t+h)-X_n(t)|_{H^{-1}} |Y_n(t+h)|_{H^1}  \nonumber  \\
      & \qquad \qquad \qquad \qquad \qquad \qquad \ \ \ + |Y_n(t+h)-Y_n(t)|_{H^{-1}} |X_n(t)|_{H^1} \bigg) dt  \nonumber \\
  \leq& \frac{1}{\g_2 n} T^{\frac 12} + C \(\bbe \int_0^{T} \sup\limits_{0<h\leq \delta}  |X_n(t+h)- X_n(t)|_{H^{-1}}^2 dt\)^{\frac 12}  \nonumber  \\
      &  + C \(\bbe \int_0^{T} \sup\limits_{0<h\leq \delta}  |Y_n(t+h)- Y_n(t)|_{H^{-1}}^2 dt\)^{\frac 12}.
\end{align}

Let us estimate   $\bbe \int_0^{T}\sup_{0<h\leq \delta} |Y_n(t+h)- Y_n(t)|_{H^{-1}}^2 dt$
in the right hand side of \eqref{L1-h-un-esti}. Similar arguments apply  to $\bbe \int_0^{T}\sup_{0<h\leq \delta} |X_n(t+h)- X_n(t)|_{H^{-1}}^2 dt$.

By the backward equation \eqref{back-equa},
\begin{align} \label{K1-K4}
   & \bbe \int_0^{T} \sup\limits_{0<h\leq \delta}|Y_n(t+h)- Y_n(t)|_{H^{-1}}^2 dt \nonumber \\
   \leq & \bbe \int_0^{T} \sup\limits_{0<h\leq \delta} \bigg| \int_t^{t+h} i\Delta Y(s) ds \bigg|_{H^{-1}}^2 dt  \nonumber \\
     & + \bbe \int_0^{T} \sup\limits_{0<h\leq \delta} \bigg| \int_t^{t+h} ( -\lbb i h_1(X_n(s) )Y_n(s)  + \lbb i h_2(X_n(s) )\ol{Y_n}(s) ) ds \bigg|_{H^{-1}}^2  dt  \nonumber \\
     & + \bbe \int_0^{T} \sup\limits_{0<h\leq \delta} \bigg| \int_t^{t+h} (\mu   - i V_0     - i u_n(s)\cdot V )Y_n(s)ds \bigg|_{H^{-1}}^2 dt \nonumber \\
     & + \bbe \int_0^{T} \sup\limits_{0<h\leq \delta} \bigg| \int_t^{t+h} ( \g_1(X_n(s) -\bbx_1(s) ) - \sum\limits_{k=1}^N \ol{\mu_k} e_kZ_{k,n}(s) )  ds \bigg|_{H^{-1}}^2 dt  \nonumber \\
    & + \bbe \int_0^{T}  \sup\limits_{0<h\leq \delta} \bigg| \int_t^{t+h} \sum\limits_{k=1}^N Z_{k,n}(s)  d\beta_k(s) \bigg|_{H^{-1}}^2 dt
   =: \sum\limits_{j=1}^5 K_j.
\end{align}

For $K_1$, by \eqref{bdd-Yn-Wpq},
\begin{align} \label{K1}
     K_1 \leq  \bbe \int_0^{T}   \sup\limits_{0<h\leq \delta}  \( \int_t^{t+h} |Y_n(s)|_{H^1} ds\)^2dt
         \leq   \delta^2 T\ \bbe \sup\limits_{t\in[0,T+1]} |Y_n(t)|^2_{H^1}
         \leq  C \delta^2,
\end{align}
where $C$ is independent of $n$.

Similarly, by Cauchy's inequality and by \eqref{bdd-Yn},
\begin{align} \label{K2}
    K_3+ K_4 \leq& \delta^2 T (|\mu|_{\9} + |V_0|_{\9} + D_U\|V\|_{L^\9(0,T+1; L^\9)}) \bbe \sup\limits_{t\in[0,T+1]}|Y_n(t)|_2^2 \nonumber \\
             & + \delta \g_1 \bbe \int_0^{T} \int_t^{t+\delta} |X_n(s)-\bbx_1(s)|_2^2 ds dt \nonumber \\
             & + \delta \sum\limits_{k=1}^N |\mu_k||e_k|_{\9} \bbe \int_0^{T} \int_t^{t+\delta} |Z_{k,n}(s)|_2^2 ds dt \nonumber \\
         \leq& C (\delta+\delta^2),
\end{align}
where $C$ is independent of $n$.

Regarding $K_2$, choose the Strichartz pair $(p,q)=(\a+1, \frac{4(\a+1)}{d(\a-1)})$.
Since $p\in(2,\frac{2d}{d-2})$, $L^{p'}(\bbr^d)\hookrightarrow H^{-1}(\bbr^d) $, we have
\begin{align*}
   K_2
   \leq&  \bbe \int_0^{T}  \sup\limits_{0<h\leq \delta}  \( \int_t^{t+h} \bigg| -\lbb i h_1(X_n(s))Y_n(s) + \lbb i h_2(X_n(s))\ol{Y_n}(s) \bigg|_{L^{p'}} ds \)^2  dt \nonumber \\
   \leq&  \a^2   \bbe \int_0^{T} \sup\limits_{0<h\leq \delta}  \( \int_t^{t+h} \bigg|  X_n^{\a-1}(s) Y_n(s)  \bigg|_{L^{p'}} ds \)^2  dt \nonumber \\
   \leq&  \delta^{2/q} \a^2  \bbe \int_0^{T}  \|X_n^{\a-1}Y_n\|^2_{L^{q'}(t,t+\delta; L^{p'})} dt \nonumber \\
   \leq&  \delta^{2/q} \a^2  T\ \bbe \| X_n^{\a-1}Y_n\|^2_{L^{q'}(0,T+1; L^{p'})}.
\end{align*}
Note that, by H\"older's inequality,
\begin{align*}
    \|X_n^{\a-1} Y_n\|_{L^{q'}(0,T; L^{p'})}
    \leq& T^\theta \|X_n\|^{\a-1}_{L^q(0,T; L^p)} \|Y_n\|_{L^q(0,T; L^p)},
\end{align*}
where $\theta =  1- d(\a-1) /4 \in (0,1)$. Hence
\begin{align*}
  K_2 \leq&  \delta^{2/q} \a^2  T^{2\theta +1}\ \bbe \| X_n\|^{2(\a-1)}_{L^q(0,T; L^p)} \|Y_n\|^2_{L^{q}(0,T+1; L^{p})}   \nonumber \\
   \leq&  \delta^{2/q} \a^2 T^{2\theta +1} \|X_n\|^{2(\a-1)}_{L^\9(\Omega; L^q(0,T; L^p))} \|Y_n\|^2_{L^2(\Omega; L^q(0,T+1; L^p))}
\end{align*}
Then, by \eqref{Cons-bdd-Xn-Lpq} and \eqref{bdd-Yn}
we obtain
\begin{align}
     K_4\leq C \delta^{2/q},
\end{align}
where $C$ is independent of $n$.

For $K_5$, using the Burkholder-Davis-Gundy inequality we get
\begin{align*}
     K_5
     \leq  C \int_0^{T} \bbe  \int_t^{t+\delta} \sum\limits_{k=1}^N |Z_{k,n}(s)|_2^2 ds dt.
\end{align*}
Then, using Fubini's theorem to interchange the sum and integrals, by \eqref{bdd-Zn} we have,
\begin{align} \label{K4}
    K_5
    \leq& C  \sum\limits_{k=1}^N \(\int_0^\delta\int_0^s + \int_\delta^T \int_{s-\delta}^s + \int_T^{T+\delta}\int_{s-\delta}^T \) |Z_{k,n}(s)|_2^2\ dt ds\nonumber \\
    \leq& 3 \delta C\ \sum\limits_{k=1}^N \bbe \int_{0}^{T+1} |Z_{k,n}(s)|_2^2 ds \leq C \delta,
\end{align}
where $C$ is independent of $n$.

Plugging \eqref{K1}-\eqref{K4} into \eqref{K1-K4}, since $2/q<1$ and $\delta <1$, we obtain
\begin{align*}
     \bbe \int_0^{T} \sup\limits_{0<h\leq \delta} |Y_n(t+h)- Y_n(t)|_{H^{-1}}^2 dt
     \leq C (\delta+ \delta^2 +\delta^{\frac 2 q}) \leq C \delta^{\frac 2q},
\end{align*}
where $C$ is independent of $n$.

The term $\bbe \int_0^{T} \sup_{0<h\leq \delta} |X_n(t+h)- X_n(t)|_{H^{-1}}^2 dt$ in the
right hand side of \eqref{L1-h-un-esti} can be estimated similarly.

Therefore, in view of \eqref{L1-h-un-esti} we obtain \eqref{L1-h-un} with $b = 1/q$, thereby proving the tightness of
$\{\mu_n\}$ and yielding the equality in \eqref{e2.10}.

Finally, the stochastic maximal principle \eqref{e5.1}
follows from Proposition  \ref{l4.1},  taking into account that (see \eqref{e4.11}) for the optimal $u^*$,
$\eta(u^*) + \zeta^* = 0$,
where $\zeta^* \in \mathcal{N}_{\calu_{ad^*}}(u^*)$. The proof is complete.  \hfill $\square$

\paragraph{An example.} We consider the case $m=1$ and $U=[0,\ell]$, where $\ell>0.$ Then,  equation \eqref{e5.1} reduces to

$$u^*(t)=\left\{\barr{lll}
0&\mbox{\ \ if}&{\rm Im}\int_{\rr}V(\xi)X^*(t,\xi)\ol{Y^*}(t,\xi)d\xi\le0,\vsp
\frac\ell\g_2&\mbox{\ \ if}&{\rm Im}\int_{\rr}V(\xi)X^*(t,\xi) \ol{Y^*}(t,\xi)d\xi\ge\ell,\vsp
\multicolumn{3}{l}{\frac1{\g_2}\, {\rm Im}\int_{\rr}V(\xi)X^*(t,\xi)\ol{Y^*} (t,\xi)d\xi\ \ \ \hfill\mbox{otherwise.}} \earr\right.$$
For the numerical computation of the optimal controller $u^*$, one can use the standard gradient descent algorithm  suggested by \eqref{e5.1}. Namely,
\begin{equation}\label{e5.2}
u_{n+1}   =  P_U ( \frac{1}{1+ 2\g_2 \rho_n} u_n +  \frac{2\rho_n}{1+ 2\g_2 \rho_n}{\rm Im} \int_{\rr^d}V(\xi)X_n(t,\xi) \ol{Y}_n(t,\xi)d\xi),
\end{equation}
where $\rho_n>0$ are suitable chosen and $X_n$, $Y_n$ are solutions to the forward--backward system \eqref{equa-x}, \eqref{back-equa} with $u=u_n$. \\

{\it \bf  Proof of Theorem \ref{Thm-deter}. }
The proof follows the lines as that of Theorem \ref{t2.6}.
As a matter of fact, in the deterministic case,
the analysis of the equation of variation and of the backward equation is much easier.

Similarly to Proposition \ref{l4.1}, we have
\begin{align*}
   \sup\limits_{u, v\in \calu_{ad}} (\|\vf^{u,\wt{u}}\|_{C([0,T]; H^1)} + \|\vf^{u,\wt{u}}\|_{L^q(0,T; W^{1,p})}) <\9,
\end{align*}
where $\wt{u}=v-u$, $u,v \in \calu_{ad}$, and $\vf^{u,\wt{u}}$ is the solution to the deterministic equation of variation (i.e.  \eqref{forward-equa} without $W$).
Moreover,
\begin{align} \label{deter-asy-X-y}
    \lim\limits_{\ve\to 0} \sup\limits_{t\in [0,T]}
    |\ve^{-1} (X^{u_\ve}(t)-X^u(t)) - \vf^{u,\wt{u}}(t)|_2^2 \to 0,
\end{align}
where $X^{u_\ve}$ and $X^u$ are the solutions to \eqref{equa-x} corresponding to $u_\ve:= u+\ve \wt{u}$ and $u$ respectively.

Regarding the backward deterministic equation, we can now use the reversing time arguments and the Strichartz estimates
to obtain directly the estimate below,
\begin{align*}
    \sup\limits_{u\in\calu_{ad}} ( \|Y^u\|_{L^\9(0,T; H^1)} + \|Y^u\|_{L^q(0,T; W^{1,p})} ) < \9.
\end{align*}
Based on these, one has also the directional derivative of $\Phi$ as in Proposition \ref{l4.1}, and  similarly to \eqref{L1-h-un},
the estimate below for the minimizing sequence of controls $\{u_n\}$ from Ekeland's principle,
\begin{align*}
     \limsup\limits_{n\to \9} \sup\limits_{0<h\leq \delta}
      \int_0^{T-h} |u_n(t+h) - u_n(t)|_m dt \leq C \delta^{1/q},
\end{align*}
which by the Riesz-Kolmogorov theorem implies that $\{u_n\}$ is relative compact in $L^1(0,T; \bbr^m)$,
thereby yielding the result. \hfill $\square$

\section{Appendix} \label{APPDIX}

\begin{lemma} (\cite[Lemma 2.17]{S15}) \label{Lem-Bound}
Let $T>0$ and $f\in C([0,T]; \bbr_+)$, such that
\begin{align*}
      f\leq a + b f^\a,\ \ on\ [0,T],
\end{align*}
where $a,b>0$, $\a>1$, $a<(1-\frac 1 \a) (\a b)^{-\frac{1}{\a-1}}$, and $f(0) \leq (\a b)^{-\frac{1}{\a-1}}$. Then,
\begin{align*}
      f\leq \frac{\a}{\a-1} a,\ \ on\ [0,T].
\end{align*}
\end{lemma}

{\it \bf Proof of \eqref{bdd-Xn-H1-Wpq} and \eqref{Cons-bdd-Xn-Lpq}.} For simplicity, we omit the dependence of $u$ in $X^u$.
We may assume $T\geq 1$ without loss of generality. Set
\begin{align}
    H(X(t)) : = \frac 12 |\na X(t)|_2^2 - \frac{\lbb}{\a +1} |X(t)|_{L^{\a+1}}^{\a+1}
\end{align}
As in the proof of \cite[Theorem 3.1]{2} we have for $t\in [0,T]$,
\begin{align} \label{Ito-Hami}
     &H(X(t)) - H(x) \nonumber \\
    =& - \int_0^t \( Im  \<(\na V_0  + u(s) \cdot \na V) X(s), \na X(s)\>\) ds  \nonumber \\
     &  +\int_0^t \( Re \<-\nabla(\mu X(s)),\nabla X(s) \>_2ds
       + \frac{1}{2} \sum\limits_{j=1}^N |\nabla (X(s)\phi_j)|_2^2 \) ds      \nonumber     \\
    & -\frac{1}{2}\lambda (\alpha-1)\sum\limits_{j=1}^N \int_0^t \int (Re\phi_j)^2  |X(s)|^{\alpha+1} d\xi ds
      + M(t),
\end{align}
where $M(t):= \sum^N_{j=1} \int_0^t (Re \<\nabla(\phi_j
    X(s)),\nabla X(s) \>_2
     - \lbb \int Re\phi_j |X(s)|^{\alpha+1} d\xi)d\beta_j(s)$.

Below we shall treat  the focusing and defocusing cases respectively.

$(i)$ (The focusing case $\lbb=1$.) Note that, by \cite[Lemma 3.5]{2},
\begin{align} \label{spli-lp}
    |X(t)|_{L^{\a+1}}^{\a+1}
    \leq C_\ve |X(t)|_2^p + \ve |\na X(t)|_2^2,
\end{align}
where $p = 2 \frac{2(\a+1)-d(\a-1)}{4-d(\a-1)} >2$.
As in the proof of \cite[Theorem 3.7]{2}, the first three terms on the right hand side of \eqref{Ito-Hami} are
bounded by $ C\int_0^t (|X(s)|_2^p + |X(s)|_2^2 + |\na X(s)|_2^2) ds$, where $C$ is independent of $u$. Thus, taking $\ve < \frac{\a+1}{4}$ yields
\begin{align}\label{esti-foc-hami}
    |\na X(t)|_2^2
    \leq 4H(x) + C (|X(t)|_2^p + D(t) ) + 4 M(t),
\end{align}
where $D(t) := \int_0^t ( |X(s)|_2^p + |X(s)|_2^2 +|\na X(s)|_2^2) ds$.
It follows that  for any  $\rho \geq 4$,
\begin{align} \label{naX-D-M}
     |\na X(t)|_2^{2\rho}
     \leq& C + C( |X(t)|_2^{\rho p} + D^\rho(t) + |M(t)|^\rho)
\end{align}
with $C$ independent of $u$.

Note that, by Jensen's inequality and the conservation $|X(t)|_2=|x|_2$, $t\in [0,T]$,
\begin{align} \label{D}
    \bbe \sup\limits_{s\in[0,t]} D^\rho(s)
    \leq& \bbe \sup\limits_{s\in[0,t]} s^{\rho -1} \int_0^s ( |X(r)|_2^{p\rho} + |X(r)|_2^{2\rho} +|\na X(r)|_2^{2\rho}) dr \nonumber \\
    \leq&  C(\rho, T)\(1 +  \int_0^t \bbe \sup\limits_{r\in[0,s]} |\na X(r)|_2^{2\rho}ds \).
\end{align}
Moreover, by the BDG inequality we get
\begin{align} \label{esti-M}
   & \bbe \sup\limits_{s\in[0,t]} |M(s)|^\rho  \nonumber \\
   \leq& C(\rho) \bbe \(\int_0^t \sum^N_{j=1} \(|Re \<\nabla(\phi_j
    X(s)),\nabla X(s) \>_2 |^2
        + \bigg|\int Re\phi_j |X(s)|^{\alpha+1} d\xi \bigg|^2\) ds \)^{\frac \rho 2}  \nonumber \\
   \leq& C(\rho) \bbe \(\int_0^t |\na X(s)|_2^4 + |X(s)|_2^4 + |X(s)|_{L^{\a+1}}^{2(\a+1)} ds \)^{\frac \rho 2} \nonumber  \\
   \leq& C(\rho, T) \bbe  \int_0^t |\na X(s)|_2^{2\rho} + |X(s)|_2^{2\rho } + |X(s)|_{L^{\a+1}}^{(\a+1)\rho} ds,
\end{align}
Then, using  the conservation  and \eqref{spli-lp}  one obtains the estimate
\begin{align} \label{m-lp}
   \bbe \sup\limits_{s\in[0,t]} |M(s)|^\rho \leq  C(\rho,T) \(1 + \int_0^t \bbe  \sup\limits_{r\in[0,s]} |\na X(r)|_2^{2\rho}ds\).
\end{align}

Thus, plugging \eqref{D} and \eqref{m-lp} into \eqref{naX-D-M} and using  the conservation yields
\begin{align*}
    \bbe  \sup\limits_{s\in[0,t]} |\na X(s)|_2^{2\rho}
    \leq C+ C \int_0^t  \sup\limits_{r\in[0,s]} |\na X(r)|_2^{2\rho}ds,
\end{align*}
which  implies \eqref{bdd-Xn-H1-Wpq} by Gronwall's inequality. \\

$(ii)$ (The defocusing case $\lbb=-1$.)
Similarly to \eqref{esti-foc-hami}, we have by \eqref{Ito-Hami},
\begin{align*}
   & \frac 12 |\na X(t)|_2^2 + \frac {1}{\a+1} |X(t)|_{L^{\a+1}}^{\a+1} \\
    \leq& H(x)
      + C \int_0^t (|X(s)|_2^2 +  |\na X(s)|_2^2 + |X(s)|^{\a+1}_{L^{\a+1}}) ds
      + M(t).
\end{align*}
Using the conservation and \eqref{esti-M} we get for $\rho \geq 4$,
\begin{align*}
     &\bbe \sup\limits_{s\in[0,t]} (|\na X(t)|_2^{2\rho} + |X(t)|^{(\a+1)\rho}_{L^{\a+1}}) \\
     \leq& C + C\int_0^t  \bbe \sup\limits_{r\in[0,s]} (|\na X(r)|_2^{2\rho} + |X(r)|^{(\a+1)\rho}_{L^{\a+1}}) ds,
\end{align*}
and so \eqref{bdd-Xn-H1-Wpq} follows.\\

It remains to prove \eqref{Cons-bdd-Xn-Lpq}. Indeed, in the case that $e_k$ are constants, $1\leq k\leq N$, by the rescaling
transformation $y= e^{-W}X$, we have
\begin{align} \label{equa-y}
    \partial_t y = -i\Delta y - \lbb i |y|^{\a-1} y + f(u)y,
\end{align}
where $f(u):= -i(V_0 + u\cdot V)$. Note that, the Strichartz coefficient $C_T$ is now identically a deterministic constant.
Then, arguing as in \eqref{esti-yn*-lpq} we obtain that
$\sup_{u\in\calu_{ad}}\|y^u\|_{L^\9(\Omega; L^q(0,T; L^p))}<\9$
for any Strichartz pair $(p,q)$.

As regards the estimate for $\|X^u\|_{L^\rho(\Omega; L^q(0,T; W^{1,p}))}$,  it suffices to prove that for any $\rho \geq 1$,
\begin{align} \label{esti-y-Wpq}
     \sup\limits_{u\in \calu_{ad}} \bbe \|y^u\|^\rho _{L^q(0,T; W^{1,p})} < \9,
\end{align}
where $(p,q)=(\a+1, \frac{4(\a+1)}{d(\a-1)})$.

Since $|\na (|y|^{\a-1}y)| \leq \a |y|^{\a-1} |\na y|$,
the H\"older inequality implies that
\begin{align} \label{h1-stri-nonl}
      \| |y|^{\a-1}y\|_{L^{q'}(0,t; W^{1,p'})}
   \leq 2\a t^{\theta} \|y\|^\a_{L^q(0,t; W^{1,p})},
\end{align}
where $\theta=1- \frac{d(\a-1)}{4}>0$.
Moreover,
\begin{align*}
    \|f(u)y\|_{L^1(0,t; H^1)}
    \leq  T\|f(u)\|_{L^\9(0,T; W^{1,\9})} \|y\|_{C([0,t]; H^1)},
\end{align*}

Thus, applying Strichartz estimates to \eqref{equa-y} and using the estimates above, we get
\begin{align} \label{y-h1-t}
      &\|y\|_{L^q(0,t; W^{1,p})} \nonumber  \\
    \leq& C (|x|_{H^1} + 2\a t^\theta \|y\|^\a_{L^q(0,t; W^{1,p})}
          + T \|f(u)\|_{L^\9(0,T; W^{1,\9})} \|y\|_{C([0,t];H^1)}) \nonumber  \\
    \leq& D(T) (\|y\|_{C([0,T]; H^1)} + t^\theta \|y\|^\a_{L^q(0,T; W^{1,p})} ),
\end{align}
where $D(T) = C(1+ 2\a + T \sup_{u\in \calu_{ad}}\|f(u)\|_{L^\9(\Omega; L^\9(0,T; W^{1,\9}))})$.

Then, similarly to \eqref{esti-yn*-lpq} we have
\begin{align*}
     \|y^u\|_{L^q(0,T; W^{1,p})}
     \leq& \([\frac{T}{t}]+1\)^{\frac 1 q} \frac{ \a}{\a-1} D(T) \|y^u\|_{C([0,T];H^1)},
\end{align*}
where $t= \a^{-\frac{\a}{\theta}} (\a-1)^{\frac{\a-1}{\theta}} (  \|y^u\|_{C([0,T]; H^1)} +1)^{-\frac{\a-1}{\theta}} D(T)^{-\frac{\a}{\theta}} (\leq T)$.

Therefore, taking into account \eqref{bdd-Xn-H1-Wpq}  we obtain \eqref{esti-y-Wpq}, thereby completing the proof. \hfill $\square$

\begin{lemma} \label{Lem-tight-L1}
Let $\mu_n$, $n\geq 1$, be a family of probability measures on $L^1(0,T; \bbr^m)$. Assume that
\begin{align} \label{tight-space}
  \lim\limits_{R\to \9} \limsup\limits_{n \to \9} \mu_n\bigg\{ v\in L^1(0,T; \bbr^m): \int_0^T  |v(t) |_m dt > R \bigg\} =0,
\end{align}
and for any $\ve >0$,
\begin{align} \label{tight-time}
    \lim\limits_{\delta \to 0} \limsup\limits_{n\to \9} \mu_n \bigg\{v\in L^1(0,T; \bbr^m): \sup\limits_{0<h\leq \delta}
       \int_0^{T-h} |v(t+h)-v(t)|_m dt  >\ve \bigg\} =0.
\end{align}
Then, $\{\mu_n\}_{n\geq 1}$ is tight in $L^1(0,T; \bbr^m)$.
\end{lemma}

{\it \bf Proof.}
Set $K_1(R) = \{ v\in L^1(0,T; \bbr^m): \int_0^T  |v(t) |_m dt \leq R \}$ and
$K_2(\delta, \ve) =  \{v\in L^1(0,T; \bbr^m): \sup_{0<h\leq \delta} \int_0^{T-h} |v(t+h)-v(t)|_m dt  \leq \ve  \}$,
where $R, \delta, \ve >0$.

Fix $\ve >0$. By \eqref{tight-space} there exists $N(=N(\ve))$, $R_1(=R_1(\ve))\geq 1$, such that
$\sup_{n\geq N} \mu_n ( K^c_1(R_1) ) \leq \frac{\ve}{2}$.
Since for each $n\geq 1$, $    \lim_{R\to \9} \mu_n ( K^c_1(R) ) =0, $
we can choose $R_2(=R_2 (\ve)) $ sufficiently large, such that $\sup_{1\leq n\leq N} \mu_n ( K^c_1(R_2) ) \leq \frac{\ve}{2N}$.
Thus, letting $R= R_1 \vee R_2$ we get $\sup _{n\geq 1} \mu_n ( K^c_1(R) ) \leq \ve.$

Similarly, since for each $k, n\geq 1$, $\lim _{\delta \to 0} \mu_n ( K^c_2(\delta, \frac 1k) ) =0$,
by \eqref{tight-time} and similar arguments as above, we can choose  $\delta_k >0$ sufficiently small such that
$\sup _{n\geq 1} \mu_n ( K^c_2(\delta_k, \frac 1k) ) \leq \frac{\ve}{2^{k}}$.

Then, set $K:= K_1(R) \cap (\bigcap_{k\geq 1} K_2(\delta_k, \frac 1k))$.
It follows from \cite[Theorem 1]{S87} that $K$ is relatively compact in $L^1(0,T; \bbr^m)$, and by the estimates above we have
$\sup_{n\geq 1} \mu_n(K^c) \leq 2 \ve$, which implies  the tightness of $\{\mu_n\}_{n\geq 1}$ in $L^1(0,T; \bbr^m)$. \hfill $\square$ \\

{\it \bf Acknowledgements.}  Financial support through SFB701 at Bielefeld University is gratefully acknowledged.
V. Barbu was partially supported by a CNCS UEFISCDI (Romania) grant, project DN-II-ID-DCE-2012-4-0156.
D. Zhang  is partially supported by NSFC (No. 11501362), China Postdoctoral Science Foundation
funded project (2015M581598).

\end{document}